\documentclass[leqno, a4paper, 12pt]{article} 
\usepackage{amsmath, amsfonts,  amssymb}
\usepackage[T1]{fontenc}
\usepackage{amsfonts}

\newtheorem{corollary}{Corollary}[section]

\newtheorem{example}{Example}[section]

\newtheorem{lemma}{Lemma}[section]

\newtheorem{remark}{Remark}[section]

\newtheorem{theorem}{Theorem}[section]
\numberwithin{equation}{section}
\date{}
\title{Uniform convexity of paranormed generalizations of $L^{p}$ spaces}
\author{ Justyna Jarczyk and Janusz Matkowski}
\renewcommand*{\abstract}{\scriptsize {\bf Abstract.}}

\begin{document}
\maketitle
\begin{abstract} 
For a measure space $(\Omega ,\Sigma ,\mu )$ and a bijective increasing function 
$\varphi :\left[ 0,\infty \right) \rightarrow \left[0,\infty \right)$ the $L^{p}$-like paranormed ($F$-normed) function space with the paranorm of the form $\mathbf{p}_{\varphi }(x)=\varphi ^{-1}\left(\int_{\Omega }\varphi \circ \left\vert x\right\vert d\mu \right)$   is considered. Main results give general conditions under which this space is uniformly convex. The Clarkson theorem on the uniform convexity of $L^{p}$-space is generalized. Under some specific assumptions imposed on $\varphi$ we give not only a proof of the uniform convexity but also show the formula of a modulus of convexity. We establish the uniform convexity of all finite-dimensional paranormed spaces, generated by a strictly convex bijection $\varphi$ of $[0, \infty)$. However, the {\it a contrario} proof of this fact provides no information on a  modulus of convexity of these spaces. In some cases it can be done, even an exact formula of a modulus can be proved. We show how to make it in the case when $S={\mathbb R}^2$ and $\varphi$ is given by $\varphi(t)={\rm e}^t-1$.
\vspace{0,2cm}

\begin{flushleft}
\vspace{0,1cm}
%\end{flushleft}
%\begin{flushleft}
\end{flushleft}
\end{abstract}

\section*{Introduction}

Given a measure space $\left(\Omega ,\Sigma ,\mu \right)$ denote by $S=S\left(\Omega,\Sigma ,\mu \right)$ 
the linear real space of all $\mu$-integrable simple
functions $x:\Omega \rightarrow \mathbb{R}$. For a bijective function 
$\varphi :\left[ 0,\infty \right) \rightarrow \left[ 0,\infty \right)$ such
that $\varphi \left(0\right) =0$ define $\mathbf{p}_{\varphi }:S\rightarrow
(0,\infty )$ by 
\begin{equation*}
\mathbf{p}_{\varphi }(x)=\varphi ^{-1}\left(\int_{\Omega }\varphi \circ
\left\vert x\right\vert d\mu \right) .
\end{equation*}%
If $\varphi$ is defined by $\varphi (t)=\varphi (1)t^{p}$ with a $p\geq 1$, the functional 
$\mathbf{p}_{\varphi }$ becomes the $L^{p}$-norm in $S.$ Under weak regularity
assumptions on $\varphi$ the homogeneity condition $\mathbf{p}_{\varphi }(tx)=
t\mathbf{p}_{\varphi }(x)$, assumed for all $x\in S$ and $t>0$, forces $\varphi$ to be a power
function (\cite{JM89}, cf. also Wnuk \cite{Wnuk}). This fact implies, in
particular, that the Orlicz space cannot be normalized like the $L^{p}$ space.
To answer the question when there are non-power functions $\varphi$ such
that $\mathbf{p}_{\varphi }$ satisfies the triangle inequality in $S$ (or
when $\mathbf{p}_{\varphi }$ is a paranorm in $S$) recall the following
converse of the Minkowski inequality theorem (\cite{JMPAMS}; for more general results cf. 
\cite{JM-Indag}, \cite{JM-JMAA2008}). \\

\textit{If there are $A,B\in\Sigma$ such that 
$0<\mu \left(A\right) <1<\mu \left(B\right)\,<\infty$ and 
$\mathbf{p}_{\varphi }\left(x+y\right) \leq \mathbf{p}_{\varphi }\left(x\right) +\mathbf{p}_{\varphi }\left(y\right)$ 
for all $x,y\in S,$ then $\varphi (t)=\varphi (1)t^{p}$ for every $t\geq 0$ and some $p\geq 1$ }.

Moreover, if the underlying measure
space $(\Omega ,\Sigma ,\mu )$ satisfies one of the following two conditions:

(i) \textit{for every  $A\in \Sigma$ we have  $\mu \left(A\right) \leq 1$  or $\mu \left(A\right) =\infty$} \\(occuring,
for instance, if $\mu \left(\Omega \right) \leq 1)$,

(ii) \textit{for every  $A\in \Sigma$ we have  $\mu \left(A\right) =0$  or  $\mu \left(A\right) \geq 1$} \\
(occuring, for instance, if $\mu$ is the counting measure), then there are broad classes of non-power functions $\varphi$ such that 
$\left(S,\mathbf{p}_{\varphi }\right)$ is a paranormed space ($F$-space).
The completion $\left(\mathcal{S}^{\varphi },\mathbf{p}_{\varphi }\right)$
of the paranormed space $\left(S,\mathbf{p}_{\varphi }\right)$ is a
natural generalization of the $L^{p}$ space. It is well-known that in the
case when $p\in \left(1,\infty \right)$ the $L^{p}$ space is uniformly
convex (\cite{A.Clarc}).

In the present paper we solve the problem, proposed twenty years ago in 
\cite{JM.Graz}, to establish general conditions on $\varphi$ which ensure that,
in each of these two cases, the paranormed space 
$\left(S,\mathbf{p}_{\varphi }\right)$ (as well as $\left(\mathcal{S}^{\varphi },\mathbf{p}_{\varphi }\right)$) is uniformly convex.

Section 1 provides basic definitions and recalls some facts concerning paranorms. We remark also that each uniformly convex normed space $X$ has a nice geometrical property. Namely, it is a bead space in the sense of Pasicki \cite{Pasicki1} (cf. also \cite{Pasicki2}-\cite{Pasicki3}).

In section 2, assuming that $\mathbf{p}_{\varphi }$ is a paranorm in $S$ and
either $\varphi$ is \textit{superquadratic, }or a related two variable
function satisfies a nice condition, we prove that $\left(S,\mathbf{p}_{\varphi }\right)$
is uniformly convex and a modulus of the convexity of this space is given
(Theorems \ref{thm1} and \ref{thm4}). These results, extending the classical Clarkson theorem
concerning the case when $\varphi$ is a power function, may be useful, for
instance, in applications of Browder - G\"{o}hde - Kirk fixed point theorem
(cf., for instance, Granas and Dugundji \cite{A.Gran}, p. 52), and its
generalizations (cf., for instance, Pasicki \cite{Pasicki1}, \cite{Pasicki2}).
Unfortunately, Theorems \ref{thm1} and \ref{thm4} do not answer the question on the uniform convexity of the space  
$\left({\mathbb R}^2,\mathbf{p}_{\varphi }\right)$, where $\varphi$ is given by $\varphi(t)={\rm e}^t-1.$ This is a basic example of the space $\left(S,\mathbf{p}_{\varphi }\right)$ generated by a measure space of type (ii) and a non-power function $\varphi$. In Theorem \ref{thm5} we prove that this space is uniformly convex giving an independent argument. 

\section{Preliminaries and auxiliary results}

Let $X$ be a real linear space. A function $\mathbf{p}:X\rightarrow \mathbb{R}$ is called a \textit{paranorm} (a \textit{total} \textit{paranorm},
Wilansky \cite{Wilansky}, p. 52; or $F$-\textit{norm, } Musielak \cite{Musielak}, p. 62) in $X$ if the following conditions are satisfied:
$\mathbf{p}\left(x\right) =0$ iff $x=0$, $\mathbf{p}\left(-x\right) =\mathbf{p}\left(x\right)$ for all $x\in X$, $\mathbf{p}$ is
subadditive, i.e.
\begin{equation*}
\mathbf{p}\left(x+y\right) \leq \mathbf{p}\left(x\right) +\mathbf{p}\left(y\right), \qquad x,y\in X,
\end{equation*}
and, if $t_{n},t\in \mathbb{R}$, $x_{n},x\in X$ for $n\in \mathbb{N}$ are
such that $t_{n}\rightarrow t,$ $\mathbf{p}\left(x_{n}-x\right) \rightarrow0$, 
then $\mathbf{p}\left(t_{n}x_{n}-tx\right) \rightarrow 0 $.

If $\mathbf{p} $ is a paranorm in $X $, then $\left(X,\mathbf{p}\right) $ is
called a \textit{paranormed space} (or  $F $\textit{-space}).

We say that a paranormed space $\left(X,\mathbf{p}\right) $ is \textit{ uniformly convex} 
if for all $r>0 $ and $\varepsilon \in \left(0,2r\right) $ there exists a 
 $\delta \left(r,\varepsilon \right) \in \left(0,r\right) $ such that
\begin{equation*}
 \mathbf{p}\left(x\right) \leq r, \,\, \mathbf{p}\left(y\right) \leq r \, \text{ and } \, \mathbf{p}\left(x-y\right)
\geq \varepsilon \quad \text{imply} \quad \mathbf{p}\left(\frac{x+y}{2}\right) \leq r-\delta \left(r,\varepsilon \right)
\end{equation*}
for all $x,y \in X $. The function $\delta :\Delta \rightarrow \left(0,\infty \right),$ where 
 $\Delta =\{ \left(r,\varepsilon \right)\in {(0, \infty)}^2 :\,  \varepsilon <2r\}$, 
is referred to as a \textit{modulus of convexity} of the space $\left(X,\mathbf{p}\right) $.

\begin{remark}\label{rem10}{\rm
It follows from Pasicki {\rm \cite{Pasicki1}, Corollary 3}, that each uniformly convex normed space 
$\left(X,\left\|\, \right\|\right)$ is a bead space, i.e. it satisfies the following nice geometrical condition: for every $r>0$, $\beta >0 $ there is a $\delta >0$ such that for every $x,y\in X$ with $\left\|x-y\right\|>\beta$ there exists a $z\in X$ such that $B\left(x,r+\delta \right) \cap B\left(y,r+\delta \right)\subset B\left(z,r-\delta \right)$. (Here $B\left(x,r\right)$ denotes the open ball centered in $x$ and with the radius $r. $)
}\end{remark}

\begin{remark}\label{rem9}
Let $(X,\mathbf{p})$ be a paranormed space and $\psi :[0,\infty)\rightarrow [0,\infty )$ be a strictly increasing and subadditive bijection. 
Then $\psi \circ \mathbf{p}$ is a paranorm in $X$. Moreover, if $(X,\mathbf{p})$ is uniformly convex with a modulus of convexity 
$\delta: \Delta\rightarrow (0, \infty)$, then $(X,\psi\circ \mathbf{p})$ (as well as its completion) is uniformly convex with the modulus of convexity $\delta_{\psi }: \Delta\rightarrow (0, \infty)$, given by
\begin{equation*}
\delta _{\psi }\left(r,\varepsilon \right) 
=r-\psi \left[\psi ^{-1}\left(r\right)-\delta \left(\psi^{-1}\left(r\right),\psi ^{-1}\left(\varepsilon \right)\right)\right].
\end{equation*}
\end{remark}
Indeed, by the monotonicity of $\psi$ and subadditivity of $\psi$ and $\mathbf{p}$ we have
\begin{eqnarray*}
\left(\psi \circ \mathbf{p}\right)\left(x+y\right) 
&=& \psi \left(\mathbf{p}\left(x+y\right)\right) 
\leq \psi \left(\mathbf{p}\left(x\right)+\mathbf{p}\left(y\right) \right) 
\leq \psi\left(\mathbf{p}\left(x\right)\right) +\psi \left(\mathbf{p}\left(y\right)\right) \\
&=&\left(\psi \circ \mathbf{p}\right)\left(x\right) +\left(\psi \circ \mathbf{p}\right)\left(y\right)
\end{eqnarray*}
for all $x,y \in X$, and it is easy to verify that $\left(\psi \circ \mathbf{p}\right)$ satisfies the
remaining properties of the paranorm. Now assume that the space $(X,\mathbf{p})$ is uniformly convex with a modulus of convexity 
$\delta:\Delta \rightarrow (0, \infty)$. Take arbitrary $(r, \varepsilon) \in \Delta$ and $x,y\in X$
such that 
\begin{equation*}
\left(\psi \circ \mathbf{p}\right)\left(x\right) \leq r,\quad \left(\psi \circ \mathbf{p}\right)\left(y\right) \leq r\quad  \text{and} \quad \left(\psi \circ \mathbf{p}\right)\left(x-y\right) \geq \varepsilon,
\end{equation*}
that is $\mathbf{p}\left(x\right) \leq \psi ^{-1}\left(r\right)$, $\mathbf{p}\left(y\right) \leq \psi ^{-1}\left(r\right)$ and 
$\mathbf{p}\left(x-y\right) \geq \psi ^{-1}\left(\varepsilon \right)$.
Making use of the uniform convexity of the space $\left(X,\mathbf{p}\right)$ we get
\begin{equation*}
\mathbf{p}\left(\frac{x+y}{2}\right) 
\leq \psi ^{-1}\left(r\right)-\delta \left(\psi ^{-1}\left(r\right),\psi ^{-1}\left(\varepsilon \right) \right),
\end{equation*}
whence, by the monotonicity of $\psi$, 
\begin{eqnarray*}
\left(\psi \circ \mathbf{p}\right)\left(\frac{x+y}{2}\right) 
&=&\psi \left(\mathbf{p}\left(\frac{x+y}{2}\right)\right) 
\leq \psi \left[ \psi^{-1}\left(r\right)-\delta\left(\psi ^{-1}\left(r\right),\psi^{-1}\left(\varepsilon \right) \right) \right] \\ 
&=&r-\delta _{\psi }\left(r,\varepsilon \right).
\end{eqnarray*}
The monotonicity of $\psi $ implies that $\delta _{\psi }\left(r,\varepsilon \right) \in \left(0,r\right),$ so $\delta _{\psi }$ is a
modulus of convexity in $(X,\psi\circ\mathbf{p})$.\\

In the rest of this section we fix a measure space $(\Omega ,\Sigma ,\mu ) $ and let $S=S(\Omega ,\Sigma ,\mu)$ be the linear real space of all $\mu $-integrable simple functions $x:\Omega \rightarrow \mathbb{R}$. Let $S_{+}=S_{+}(\Omega ,\Sigma ,\mu)=\left\{ x\in S:x\geq 0\right\} $. For a bijection  $\varphi :[0,\infty)\rightarrow [0,\infty )$ such that $\varphi \left(0\right) =0 $ the functional 
$\mathbf{p}_{\varphi }:S\rightarrow (0,\infty ) $, given by
\begin{equation*}
\mathbf{p}_{\varphi }(x)=\varphi ^{-1}\left(\int_{\Omega }\varphi \circ\left\vert x\right\vert d\mu \right),
\end{equation*}
is correctly defined. Moreover, if $x\in S $ then there exist $k\in \mathbb{N} $, $r_{i}\in \mathbb{R} $, and the pairwise disjoined sets $A_{i}\in \Sigma $, $i=1,...,k $, such that $x=\sum_{i=1}^{k}r_{i}\chi _{A_{i}} $\ and 
\begin{equation*}
\mathbf{p}_{\varphi }(x)=
\varphi ^{-1}\left(\sum_{i=1}^{k}\varphi \left(\left\vert r_{i}\right\vert \right) \mu \left(A_{i}\right) \right).
\end{equation*}
(Here $\chi _{A} $ denotes the characteristic function of the set $A $.)\\

We start with two simple facts dealing with the situation where $\varphi$ is an increasing function.

\begin{remark}\label{rem11}
Let $\varphi:[0, \infty)\rightarrow \mathbb R$ be an increasing and convex $[$strictly increasing and strictly convex$]$ function. Then the function $\varphi\circ\left|\cdot\right|$ is convex $[$strictly convex$]$.
\end{remark}

{\bf Proof.} Assume, for instance, that $\varphi$ is strictly increasing and strictly convex. Fix arbitrarily $r, s \in \mathbb R$, $r\not=s$, and $\lambda \in (0,1)$. If $rs<0$ then $\left|\lambda r +(1-\lambda)s\right|<\lambda \left|r\right| +(1-\lambda)\left|s\right|$; otherwise we would have $\left|\lambda r +(1-\lambda)s\right|=\left|\lambda r\right| +\left|(1-\lambda)s\right|$ which would mean that $\lambda r$ and $(1-\lambda)s$ had the same signs, so had $r$ and $s$. Thus
\begin{equation*}
\varphi\left(\left|\lambda r +(1-\lambda)s\right|\right)<\varphi\left(\lambda \left|r\right| +(1-\lambda)\left|s\right|\right)
\leq \lambda \varphi\left(\left|r\right|\right)+\left(1- \lambda\right) \varphi\left(\left|s\right|\right).
\end{equation*}
If $rs\geq 0$ then $\left|r\right|\not=\left|s\right|$, and thus
\begin{equation*}
\varphi\left(\left|\lambda r +(1-\lambda)s\right|\right)\leq\varphi\left(\lambda \left|r\right| +(1-\lambda)\left|s\right|\right)
< \lambda \varphi\left(\left|r\right|\right)+\left(1- \lambda\right) \varphi\left(\left|s\right|\right).
\end{equation*}
Consequently, $\varphi\circ\left|\cdot\right|$ is strictly convex. $\square$

\begin{corollary}\label{c12}
Let $\varphi:[0, \infty)\rightarrow [0, \infty)$ be a convex bijection. Then for any  $r \in (0, \infty)$ the set $\left\{x \in S\colon \mathbf{p}_{\varphi }(x)\leq r\right\}$ is convex. 
\end{corollary}

{\bf Proof.} 
It follows from the assumptions that $\varphi$ is increasing. Fix arbitrarily $r \in (0, \infty)$ and take any points $x, y \in S$ satisfying $\mathbf{p}_{\varphi }(x)\leq r$, $\mathbf{p}_{\varphi }(y)\leq r$, and a number $\lambda \in (0,1)$. Then, by Remark \ref{rem11},
\begin{eqnarray*}
 \varphi \left(\mathbf{p}_{\varphi }\left(\lambda x +(1-\lambda)y \right)\right)
& = & \int_{\Omega}\varphi\circ \left|\lambda x +(1-\lambda)y\right|{\rm d}\mu \\
& \leq &  \int_{\Omega}\left(\lambda \varphi\circ \left|x\right| +(1-\lambda)\varphi\circ\left|y\right|\right){\rm d}\mu \\
& = & \lambda \int_{\Omega}\varphi\circ\left|x\right|{\rm d}\mu+\left(1- \lambda\right) \int_{\Omega} \varphi\circ\left|y\right|{\rm d}\mu \\
& = &\lambda\varphi\left(\mathbf{p}_{\varphi }(x)\right)+\left(1- \lambda\right)\varphi\left(\mathbf{p}_{\varphi }(y)\right) \\
& \leq &  \lambda \varphi(r)+\left(1- \lambda \right)\varphi (r)=\varphi (r),
\end{eqnarray*}
that is $\mathbf{p}_{\varphi }\left(\lambda x +(1-\lambda)y \right)\leq r$, which was to be proved. $\square$

\vspace{0,2cm}

From \cite{JM.Graz}, Theorem 6(a) (cf. also Hardy, Littlewood, P\'{o}lya \cite{HLP}, Theorem 106(ii)), we have the following

\begin{lemma}\label{lem1}
Let  $\varphi :[0,\infty)\rightarrow [0,\infty )$ be an increasing bijection.

{\rm (a)} Assume that $\mu (\Omega )=1$ and there is a set $A\in \Sigma$ such
that $0<\mu \left(A\right) <1$. Then $\mathbf{p}_{\varphi }$ is a paranorm
in $S(\Omega ,\Sigma ,\mu )$ if and only if the function $F:[0,\infty)^{2}\rightarrow [0,\infty )$, given by
\begin{equation}\label{e1}
F\left(r,s\right) =\varphi \left(\varphi ^{-1}\left(r\right) +\varphi^{-1}\left(s\right) \right),
\end{equation}
is concave.

{\rm (b)} If $\mu (\Omega )\leq 1$ and the function $F$ is concave, then $\mathbf{p}_{\varphi }$ is a paranorm in $S(\Omega ,\Sigma ,\mu )$.
\end{lemma}

\begin{remark}\label{rem1}{\rm (\cite{JM.Graz}, proof of Theorem 13)}.  
Let $\varphi :[0,\infty )\rightarrow [0,\infty )$ be an increasing bijection. Assume that the function 
$F:[0,\infty )^2\rightarrow [0,\infty )$, defined by \eqref{e1}, is concave. Then $\varphi $ is continuous, $\varphi \left(0\right) =0$, and the functional $\mathbf{p}_{\varphi }: S\rightarrow [0,\infty )$ has the following properties: $\mathbf{p}_{\varphi }\left(x\right) =0$ for an $x\in S$ iff 
$x=0$ $\mu $-a.e.; $\mathbf{p}_{\varphi }\left(x\right) =\mathbf{p}_{\varphi}\left(-x\right)$ 
for all $x\in S; $ if $t_{n},t\in \mathbb{R}$, 
$x_{n},x\in S$ for $n\in \mathbb{N}$ are such that 
$t_{n}\rightarrow t, \mathbf{p}_{\varphi }\left(x_{n}-x\right) \rightarrow 0$, then 
$\mathbf{p}_{\varphi }\left(t_{n}x_{n}-tx\right) \rightarrow 0$.
\end{remark}

\begin{remark}\label{rem2} 
Let $\varphi :[0,\infty )\rightarrow [0,\infty )$ be an increasing bijection. Assume that $\mu (\Omega )=1$ and there is a set 
$A\in \Sigma$ such that $0<\mu \left(A\right) <1$. Then the function $F:[0,\infty )^2\rightarrow [0,\infty )$, given by \eqref{e1}, is convex iff $\mathbf{p}_{\varphi }(x+y)\geq \mathbf{p}_{\varphi }(x)+\mathbf{p}_{\varphi }(y)$ 
for all $x,y\in S_{+}\left(\Omega ,\Sigma ,\mu \right)$.
\end{remark}

\begin{remark}\label{rem3} {\rm (\cite{JMPAMS}).}  
Let $\varphi :\left[ 0,\infty \right) \rightarrow \left[0,\infty \right)$ be a twice differentiable function mapping $[0, \infty)$ onto itself, such that $\varphi ^{\prime }\left(r\right) >0$ and $\varphi^{\prime \prime }\left(r\right) \neq 0$ for all $r>0$, and let $F:\left[0,\infty \right)^2 \rightarrow \left[0,\infty \right)$  be defined by \eqref{e1}.

The function $F$ is concave if and only if $\varphi ^{\prime \prime }\left(r\right) >0$ for all $r>0$ and the function 
$\frac{\varphi ^{\prime }}{\varphi ^{\prime \prime}}$ is superadditive in $\left(0,\infty \right)$, i.e.
\begin{equation*}
\frac{\varphi ^{\prime }\left(r+s\right) }{\varphi ^{\prime \prime }\left(r+s\right)}
\geq \frac{\varphi ^{\prime }\left(r\right) }{\varphi ^{\prime \prime }\left(r\right)}
+\frac{\varphi ^{\prime }\left(s\right) }{\varphi^{\prime \prime }\left(s\right)}, \qquad r,s>0.
\end{equation*}
The function $F$ is convex if and only if $\varphi ^{\prime \prime }\left(r\right) <0$ for all $r>0$ and the
function $\frac{\varphi ^{\prime }}{\varphi ^{\prime \prime }}$ is subadditive in $\left(0,\infty \right)$, i.e. 
\begin{equation*}
\frac{\varphi ^{\prime }\left(r+s\right) }{\varphi ^{\prime \prime }\left(r+s\right) }
\leq \frac{\varphi ^{\prime }\left(r\right) }{\varphi ^{\prime\prime }\left(r\right) }
+\frac{\varphi ^{\prime }\left(s\right) }{\varphi^{\prime \prime }\left(s\right) }, \qquad r,s>0.
\end{equation*}
\end{remark}

\begin{lemma}\label{lem2}
Let $\varphi :\left[0,\infty \right) \rightarrow \left[ 0,\infty \right)$ 
be a twice differentiable function such that $\varphi \left(0\right) =0$ 
and $\varphi ^{\prime }\left(r\right) >0, $ $\varphi ^{\prime \prime }\left(r\right) >0$ 
for all $r>0 $. Let $G:\left[ 0,\infty \right) ^{2}\rightarrow \left[ 0,\infty \right)$ 
be defined by
\begin{equation*}
G\left(r,s\right)=\varphi \left(\left\vert \varphi ^{-1}\left(r\right)
-\varphi ^{-1}\left(s\right) \right\vert \right).
\end{equation*}
If the function $\frac{\varphi ^{\prime }}{\varphi ^{\prime \prime }} $ is
superadditive in $\left(0,\infty \right)$, then $G $ is convex. 
\end{lemma}

{\bf Proof.} 
Put $D_{1}:=\left\{ \left(r,s\right) \in \left[ 0,\infty \right) ^{2}:r\geq s\right\}$ 
and $D_{2}:=\{ \left(r,s\right) \in \left[ 0,\infty \right) ^{2}:r\leq s \}$. 
We first show that $G $ is convex in the interior of $D_{1} $. Note that for arbitrary $r>s>0$ and $u,v\in \mathbb{R}$ the function $g_{r,s,u,v}$, defined by
\begin{equation*}
g_{r,s,u,v}\left(t\right)=G\left(r+tu,s+tv\right) =\varphi \left(\varphi ^{-1}\left(r+tu\right) -\varphi ^{-1}\left(s+tv\right) \right),
\end{equation*}
is twice differentiable in a neighbourhood of $0$. To prove the convexity of $G$ in the interior of $D_1$ it is enough to show that 
$g_{r,s,u,v}^{\prime \prime }\left(0\right) \geq 0$ for all $r>s>0$ and $u,v \in {\mathbb R}$ (cf. \cite{HLP}, p. 86). 
We have 
\begin{eqnarray*}
&& \hspace{-0,7cm} g_{r,s,u,v}^{\prime }\left(t\right) \\ 
&&  =\varphi ^{\prime }\left(\varphi^{-1}\left(r+tu\right) -\varphi ^{-1}\left(s+tv\right) \right) 
\left[ \frac{u}{\varphi ^{\prime }\left(\varphi ^{-1}\left(r+tu\right) \right) }-
			 \frac{v}{\varphi ^{\prime }\left(\varphi ^{-1}\left(s+tv\right) \right) }\right] 
\end{eqnarray*}
and
\begin{eqnarray*}
&& \hspace{-0,7cm} g_{r,s,u,v}^{\prime \prime }\left(t\right)  \\
&& \hspace{-0,7cm}=\varphi ^{\prime \prime }\left(\varphi ^{-1}\left(r+tu\right) -\varphi ^{-1}\left(s+tv\right) \right) 
		\left[ \frac{u}{\varphi ^{\prime }\left(\varphi ^{-1}\left(r+tu\right)\right) }+
	    		 \frac{v}{\varphi ^{\prime }\left(\varphi^{-1}\left(s+tv\right)\right) }\right] ^{2}\\
&& \hspace{0cm}+ \varphi ^{\prime }\left(\varphi ^{-1}\left(r+tu\right) -\varphi^{-1}\left(s+tv\right) \right)\\ 
&&	\hspace{3,5cm}	\cdot\left[ \frac{\varphi ^{\prime \prime }\left(\varphi ^{-1}\left(s+tv\right) \right) v^{2}}
		            {\left(\varphi ^{\prime }\left(\varphi ^{-1}\left(s+tv\right) \right)\right) ^{3}}-
		       \frac{\varphi ^{\prime \prime}\left(\varphi ^{-1}\left(r+tu\right) \right) u^{2}}
		            {\left(\varphi^{\prime }\left(\varphi ^{-1}\left(r+tu\right) \right) \right) ^{3}}\right]
\end{eqnarray*}
for all $t$ from a neighbourhood of $0.$ Since $\varphi $ is an increasing bijection of $\left(0,\infty \right)$, it is enough to show that 
$g_{\varphi \left(r\right),\varphi \left(s\right),u,v}^{\prime \prime}\left(0\right) \geq 0$ for all $r>s>0$ and $u,v\in \mathbb{R}$. Putting $g=g_{\varphi \left(r\right) ,\varphi \left(s\right) ,u,v}$ we have
\begin{equation*}
g^{\prime \prime }\left(0\right) =Au^{2}+2Buv+Cv^{2},
\end{equation*}
where
\begin{eqnarray*}
&& A =\frac{\varphi ^{\prime \prime }\left(r-s\right) }{\left[ \varphi^{\prime }\left(r\right) \right] ^{2}}-
\frac{\varphi ^{\prime }\left(r-s\right) \varphi ^{\prime \prime }\left(r\right) }
{\left[ \varphi^{\prime }\left(r\right) \right]^{3}}, \qquad 
B=\frac{\varphi^{\prime \prime }\left(r-s\right) }{\varphi ^{\prime }\left(r\right)\varphi ^{\prime }\left(s\right) }, \\
&& C = \frac{\varphi ^{\prime \prime }\left(r-s\right) }{\left[\varphi ^{\prime }\left(s\right) \right] ^{2}}+
\frac{\varphi ^{\prime}\left(r-s\right) \varphi ^{\prime \prime }\left(s\right) }{\left[ \varphi^{\prime }\left(s\right) \right] ^{3}}.
\end{eqnarray*}
Since  
\begin{eqnarray*}
&& \hspace{-1cm} AC-B^{2} \\
&& \hspace{0cm} = -\frac{\left[ \varphi ^{\prime }\left(r-s\right) \right] ^{2}\varphi^{\prime \prime }\left(r\right) \varphi ^{\prime \prime }\left(s\right) }{\left[ \varphi ^{\prime }\left(r\right) \right] ^{3}\left[ \varphi ^{\prime}\left(s\right) \right] ^{3}}
+\frac{\varphi ^{\prime \prime }\left(r-s\right) \varphi ^{\prime }\left(r-s\right) \varphi ^{\prime \prime}\left(s\right) }
			{\left[ \varphi ^{\prime }\left(r\right) \right] ^{2}\left[\varphi ^{\prime }\left(s\right) \right] ^{3}}\\
&& \hspace{7,0cm} -\frac{\varphi^{\prime \prime }\left(r-s\right)\varphi ^{\prime }\left(r-s\right)\varphi ^{\prime \prime }\left(r\right) }
			{\left[ \varphi ^{\prime }\left(r\right) \right] ^{3}\left[ \varphi ^{\prime }\left(s\right) \right] ^{2}}\\
&& \hspace{0cm} =\frac{\varphi ^{\prime }\left(r-s\right) \varphi ^{\prime \prime }\left(r-s\right) 
								\varphi ^{\prime \prime }\left(r\right) \varphi ^{\prime \prime}\left(s\right) }
								{\left[ \varphi ^{\prime }\left(r\right) \right] ^{3}\left[ \varphi ^{\prime }\left(s\right) \right]^{3}}
								\left(\frac{\varphi^{\prime }\left(r\right) }{\varphi ^{\prime \prime }\left(r\right) }
								-\frac{\varphi ^{\prime }\left(r-s\right) }{\varphi ^{\prime \prime }\left(r-s\right) }
								-\frac{\varphi ^{\prime }\left(s\right) }{\varphi ^{\prime\prime }\left(s\right) }\right) ,
\end{eqnarray*}
the superadditivity of $\frac{\varphi ^{\prime }}{\varphi ^{\prime \prime }}$ implies that $AC-B^{2}\geq 0 $. The superadditivity and positivity of $\frac{\varphi ^{\prime }}{\varphi ^{\prime \prime }}$ imply that $\frac{\varphi^{\prime }}{\varphi ^{\prime \prime }}$ is increasing. Thus the function $\frac{\varphi ^{\prime \prime }}{\varphi ^{\prime }}$ is decreasing, whence 
\begin{equation*}
A=\frac{\varphi ^{\prime \prime }\left(r-s\right) }{\left[ \varphi ^{\prime}\left(r\right) \right] ^{2}}
-\frac{\varphi ^{\prime}\left(r-s\right)\varphi^{\prime \prime}\left(r\right)}{\left[\varphi^{\prime }\left(r\right)\right]^{3}}
=\frac{\varphi^{\prime }\left(r-s\right)}{\left[\varphi ^{\prime }\left(r\right)\right] ^{2}}
  \left(\frac{\varphi^{\prime\prime }\left(r-s\right) }{\varphi ^{\prime }\left(r-s\right)}
  -\frac{\varphi ^{\prime \prime }\left(r\right) }{\varphi ^{\prime }\left(r\right)}\right) \geq 0.
\end{equation*}
Analogously, we infer that $C\geq 0$. Consequently, $g^{\prime \prime }\left(0\right) \geq 0$. This proves that the function $G $ is convex in the interior of $D_{1}.$ The convexity of $G$ in $D_{1}$ follows from the continuity of $G$, whereas the same property of $G$ in $D_{2}$ is an immediate consequence of the symmetry of $G$ and its convexity in $D_{1}$.

Note that $D_{1}\cap D_{2}=\left\{ \left(r,r\right): r\geq 0\right\}$ and $G\left(r,r\right)=0$ for all $r\geq 0. $ Take any two points $\left(r_{1},s_{1}\right) \in {\rm int}D_{1}$ and $\left(r_{2},s_{2}\right)\in {\rm int}D_{2}$. The function 
$\gamma: \left[0,1\right] \rightarrow \mathbb{R}$, defined by
\begin{equation*}
\gamma \left(t\right)=G\left(t\left(r_{2},s_{2}\right)+\left(1-t\right) \left(r_{1},s_{1}\right)\right),
\end{equation*}
is nonnegative and there is a unique $t_{0}\in \left(0,1\right)$ such that 
$t_{0}\left(r_{2},s_{2}\right) +\left(1-t_{0}\right)\\ \cdot \left(r_{1},s_{1}\right) =\left(r,r\right) \in $ $D_{1}\cap D_{2}$. By the
previous part of the proof the function $\gamma$ is convex in each of the intervals $\left[ 0,t_{0}\right] $ and $\left[t_{0},1\right]$. Since $\gamma \left(t_{0}\right)=0 $, the function $\gamma $ attains its global minimum at $t_{0}$. It follows that $\gamma $ is convex in 
$\left[0,1\right]$ and, consequently, the function $G $ is convex.
$\square $

\begin{remark}\label{rem4}{\rm 
Let $\varphi : \left[0,\infty \right) \rightarrow \left[0,\infty \right)$ be a twice differentiable function such that $\varphi(0)=0$, and  
$\varphi ^{\prime }(r)>0$, $\varphi ^{\prime \prime }(r)>0$, for all $r>0$. In Corollary \ref{cor3} the assumption of the convexity of the function 
$H:\left[0,\infty \right)^{2}\rightarrow \left[ 0,\infty \right)$, defined by
\begin{equation*}
H\left(r,s\right)=\varphi \left(\varphi ^{-1}\left(r\right)+\varphi^{-1}\left(s\right)\right)
+\varphi \left(\left\vert \varphi ^{-1}\left(r\right) -\varphi ^{-1}\left(s\right) \right\vert \right),
\end{equation*}
occurs. To establish conditions on $\varphi $ under which the function $H $ is convex we can proceed similarly as in the proof of 
Lemma \ref{lem2}. For arbitrary $r>s>0$ and $u,v\in \mathbb{R}$ the function $h_{r,s,u,v}$, defined by 
\begin{equation*}
h_{r,s,u,v}\left(t\right)=H\left(r+tu,s+tv\right),
\end{equation*}
is twice differentiable in a neighbourhood of $0$. To get the convexity of $H$ it is enough to know that $h_{\varphi(r),\varphi(s),u,v}^{\prime \prime}(0)\geq 0$ for all $r>s>0$ and $u,v \in \mathbb R$. 
Calculating $h_{r,s,u,v}^{\prime \prime }$ and  putting $h=h_{\varphi(r),\varphi(s),u,v}$ we have 
\begin{equation*}
h^{\prime \prime }\left(0\right) =Au^{2}+2Buv+Cv^{2},
\end{equation*}
where 
\begin{eqnarray*}
&& A=\frac{\varphi ^{\prime \prime}\left(r+s\right)}{\left[\varphi ^{\prime}\left(r\right)\right]^{2}}
      -\frac{\varphi^{\prime }\left(r+s\right)\varphi ^{\prime \prime }\left(r\right)}
            {\left[ \varphi ^{\prime }\left(r\right) \right] ^{3}}
      +\frac{\varphi ^{\prime \prime }\left(r-s\right)}{\left[ \varphi ^{\prime }\left(r\right) \right] ^{2}}
      -\frac{\varphi^{\prime}\left(r-s\right)\varphi^{\prime \prime}\left(r\right)}{\left[\varphi ^{\prime }\left(r\right) \right] ^{3}},\\
&& B=\frac{\varphi^{\prime\prime}\left(r+s\right)-\varphi^{\prime\prime}\left(r-s\right)}
			{\varphi^{\prime}\left(r\right)\varphi^{\prime}\left(s\right)},\\
&& C=\frac{\varphi^{\prime \prime}\left(r+s\right)}{\left[\varphi ^{\prime}\left(s\right) \right] ^{2}}
      -\frac{\varphi^{\prime}\left(r+s\right)\varphi^{\prime\prime}\left(s\right)}{\left[\varphi^{\prime}\left(s\right)\right]^{3}}
      +\frac{\varphi ^{\prime \prime }\left(r-s\right)}{\left[ \varphi ^{\prime }\left(s\right) \right] ^{2}}
      +\frac{\varphi^{\prime}\left(r-s\right)\varphi^{\prime\prime}\left(s\right)}{\left[\varphi ^{\prime }\left(s\right)\right] ^{3}}.
\end{eqnarray*}
We need to know that $AC-B^{2}\geq 0, A\geq 0$ and $C\geq 0$.

Since 
\begin{eqnarray*}
&& \hspace{-0,55cm} \frac{\left(AC-B^{2}\right)\left[\varphi ^{\prime }\left(r\right)\right]^{2}\left[\varphi^{\prime }\left(s\right)\right]^{2}}       {\left[\varphi^{\prime }\left(r+s\right)\right]^2 -\left[\varphi ^{\prime }\left(r-s\right)\right]^2}
=4\frac{\varphi^{\prime \prime}\left(r+s\right)\varphi ^{\prime \prime}\left(r-s\right)}
     {\left[\varphi^{\prime }\left(r+s\right)\right]^2 -\left[\varphi ^{\prime }\left(r-s\right)\right]^{2}}
 +\frac{\varphi^{\prime \prime}\left(r\right)}{\varphi^{\prime}\left(r\right)}
      \frac{\varphi^{\prime \prime}\left(s\right)}{\varphi^{\prime }\left(s\right)}\\
&& \hspace{1cm}
-\frac{\varphi ^{\prime \prime }\left(r\right) }{\varphi ^{\prime }\left(r\right)} \cdot   
 \frac{\varphi^{\prime\prime}\left(r+s\right)+\varphi^{\prime\prime}\left(r-s\right)}
       {\varphi^{\prime}\left(r+s\right)-\varphi^{\prime}\left(r-s\right)}
-\frac{\varphi ^{\prime \prime }\left(s\right) }{\varphi^{\prime }\left(s\right)} \cdot
 \frac{\varphi^{\prime\prime}\left(r+s\right)+\varphi^{\prime\prime}\left(r-s\right)}
       {\varphi^{\prime}\left(r+s\right)+\varphi^{\prime}\left(r-s\right)},
\end{eqnarray*}
we conclude that $AC-B^{2}\geq 0 $ iff the function on the right-hand side of the above equality is nonnegative. Moreover, $A\geq 0$ iff  
\begin{equation*}
\left[\varphi ^{\prime \prime }\left(r+s\right)+\varphi ^{\prime \prime}\left(r-s\right)\right]\varphi ^{\prime }\left(r\right) \geq \left[\varphi ^{\prime }\left(r+s\right) +\varphi ^{\prime }\left(r-s\right)\right]\varphi ^{\prime \prime }\left(r\right),
\end{equation*}
and $C\geq 0 $ iff
\begin{equation*}
\left[\varphi^{\prime \prime}\left(r+s\right)+\varphi^{\prime \prime}\left(r-s\right)\right]\varphi ^{\prime }\left(s\right) 
\geq \left[\varphi ^{\prime }\left(r+s\right) -\varphi ^{\prime }\left(r-s\right)\right] \varphi ^{\prime \prime }\left(s\right).
\end{equation*}
Consequently, the following inequalities, held for all $r>s>0$, provide a sufficient condition for $H$ to be convex:
\begin{eqnarray*}
&& \hspace{-0,75cm} 4\frac{\varphi^{\prime \prime}\left(r+s\right)\varphi ^{\prime \prime}\left(r-s\right)}
     {\left[\varphi^{\prime }\left(r+s\right)\right]^2 -\left[\varphi ^{\prime }\left(r-s\right)\right]^{2}}
 +\frac{\varphi^{\prime \prime}\left(r\right)}{\varphi^{\prime}\left(r\right)}
      \frac{\varphi^{\prime \prime}\left(s\right)}{\varphi^{\prime }\left(s\right)}\\
&& \hspace{0,45cm}
\geq \frac{\varphi ^{\prime \prime }\left(r\right) }{\varphi ^{\prime }\left(r\right)} \cdot   
 \frac{\varphi^{\prime\prime}\left(r+s\right)+\varphi^{\prime\prime}\left(r-s\right)}
       {\varphi^{\prime}\left(r+s\right)-\varphi^{\prime}\left(r-s\right)}
+\frac{\varphi ^{\prime \prime }\left(s\right) }{\varphi^{\prime }\left(s\right)} \cdot
 \frac{\varphi^{\prime\prime}\left(r+s\right)+\varphi^{\prime\prime}\left(r-s\right)}
       {\varphi^{\prime}\left(r+s\right)+\varphi^{\prime}\left(r-s\right)},
\end{eqnarray*}
\begin{eqnarray*}
\frac{\varphi^{\prime\prime}\left(r+s\right)+\varphi^{\prime\prime}\left(r-s\right)}
       {\varphi^{\prime}\left(r+s\right)+\varphi^{\prime}\left(r-s\right)}
\geq \frac{\varphi ^{\prime \prime }\left(r\right) }{\varphi ^{\prime }\left(r\right)}, 
\end{eqnarray*}
and
\begin{eqnarray*}
\frac{\varphi^{\prime\prime}\left(r+s\right)+\varphi^{\prime\prime}\left(r-s\right)}
			{\varphi^{\prime}\left(r+s\right)-\varphi^{\prime}\left(r-s\right)}
 \geq \frac{\varphi ^{\prime \prime }\left(s\right) }{\varphi ^{\prime }\left(s\right)}. 
\end{eqnarray*}
We omit a simple calculation showing that if $\varphi(t)=t^p, t\geq 0$, with a $p\in (1,2)$, then the above three inequalities are satisfied, so $H$ is convex in that case. 
}\end{remark}

The next lemma contains the Mulholland inequality \cite{Mulhol} as a special case.

\begin{lemma}\label{lem3}{\rm (\cite{JM-Abh}, \cite{JM-JMAA-2012})}. 
Assume that for every $A \in \Sigma$ either $\mu(A)=0$, or $\mu(A)\geq 1$. 
If $\varphi :[0,\infty )\rightarrow [0,\infty )$ is an increasing bijection, convex, and geometrically convex, i.e. 
\begin{equation*}
\varphi \left(\sqrt{st}\right) \leq \sqrt{\varphi \left(s\right) \varphi \left(t\right)}, \qquad s,t\geq 0,
\end{equation*}
then $\mathbf{p}_{\varphi }$ is a paranorm in $S(\Omega,\Sigma ,\mu)$.
\end{lemma}

One can easily observe the following fact providing a simple tool when examining the geometric convexity of a function. It is based on an obvious fact that a function $\varphi:(0, \infty)\rightarrow(0, \infty)$ is geometrically convex if and only if the function $\log \circ \varphi \circ {\rm exp}$ is convex.

\begin{remark}\label{rem12}
A differential function $\varphi:(0, \infty)\rightarrow(0, \infty)$  is geometrically convex iff the function  
\begin{equation*}
(0, \infty)\ni t\rightarrow t\frac{\varphi'(t)}{\varphi(t)}
\end{equation*}
is increasing.
\end{remark}

\begin{remark}\label{rem5}{\rm
Making use of Remark \ref{rem12} it is easy to verify that the functions $\varphi:[0, \infty)\rightarrow [0, \infty)$, defined by $\varphi \left(t\right)=t^{p}$ with a $p\geq 1$, as well as the non-power functions of the form $\varphi \left(t\right)=a^{t}-1$ with an $a>1,$ and $\varphi (t)=t^{p}a^{t}$
with $a>1$ and $p\geq 1$  satisfy the assumptions of Lemma \ref{lem3}.
}\end{remark}

\section{Uniform convexity of the space $\left(S(\Omega ,\Sigma , \mu ),\mathbf{p}_{\varphi }\right)$}

In the whole section $(\Omega ,\Sigma ,\mu )$ means a measure space. We begin with the following

\begin{lemma}\label{lem4}
Let $\varphi :[0,\infty)\rightarrow [0,\infty )$ be a bijection with $\varphi \left(0\right)=0$. If $\varphi $ is superquadratic, i.e.
\begin{equation*}
\varphi \left(r+s\right) +\varphi \left(\left\vert r-s\right\vert \right)\geq 2\varphi \left(r\right) +2\varphi \left(s\right), 
\qquad r,s >0,
\end{equation*}
then 
\begin{equation*}
\varphi \left(\mathbf{p}_{\varphi }\left(x+y\right)\right)+\varphi\left(\mathbf{p}_{\varphi }\left(x-y\right)\right) 
\geq 2\varphi \left(\mathbf{p}_{\varphi }\left(x\right) \right) +2\varphi \left(\mathbf{p}_{\varphi }\left(y\right) \right)
\end{equation*}
for all $x,y\in S(\Omega ,\Sigma ,\mu )$. 
\end{lemma}

{\bf Proof.} Take arbitrary $x,y \in S(\Omega ,\Sigma ,\mu )$. Then 
\begin{eqnarray*}
\varphi\left(\left|x(t)+y(t)\right|\right)+\varphi\left(\left|x(t)-y(t)\right|\right)
\geq 2\varphi\left(\left|x(t)\right|\right)+2\varphi\left(\left|y(t)\right|\right), \quad t \in \Omega,
\end{eqnarray*}
whence
\begin{eqnarray*}
\int_{\Omega}\varphi\circ\left|x+y\right|{\rm d}\mu +\int_{\Omega}\varphi\circ\left|x-y\right|{\rm d}\mu
\geq 2\int_{\Omega}\varphi\circ\left|x\right|{\rm d}\mu +2\int_{\Omega}\varphi\circ\left|y\right|{\rm d}\mu, 
\end{eqnarray*}
that is
\begin{equation*}
\varphi \left(\mathbf{p}_{\varphi }\left(x+y\right) \right) +\varphi\left(\mathbf{p}_{\varphi }\left(x-y\right) \right) 
\geq 2\varphi \left(\mathbf{p}_{\varphi }\left(x\right) \right) +2\varphi \left(\mathbf{p}_{\varphi }\left(y\right) \right),
\end{equation*}
which was to be shown.
 $\square $

\begin{remark}\label{rem6}
For $p\geq 2$ the power function $\varphi: [0, \infty)\rightarrow [0, \infty)$, given by $\varphi \left(t\right) =t^{p}$, is superquadratic, and for $p\in \left(0,2\right]$ it is subquadratic, i.e. $\varphi $ satisfies the inequality
\begin{equation*}
\varphi\left(r+s\right)+\varphi\left(\left\vert r-s\right\vert\right)\leq 2\varphi\left(r\right)+2\varphi\left(s\right), \quad r,s\geq 0.
\end{equation*}
\end{remark}

{\bf Proof.} 
For $p=2 $ this is obvious. Let us fix $p>2 $. It is enough to show that $\left(x+y\right)^{p}+\left(x-y\right)^{p}\geq 2x^{p}+2y^{p}$ for all $x\geq y \geq 0$. Setting $t=\frac{y}{x} $ we can write this inequality in the equivalent form
\begin{equation*}
\left(1+t\right)^{p}+\left(1-t\right)^{p}\geq 2+2t^{p}, \quad t\in [0,1].
\end{equation*}
Put $f\left(t\right)=\left(1+t\right)^{p}+\left(1-t\right)^{p}-2-2t^{p}$ for $t\in \left[ 0,1\right]$. Then 
\begin{equation*}
f^{\prime }\left(t\right)=p\left[\left(1+t\right)^{p-1}-\left(1-t\right) ^{p-1}-2t^{p-1}\right] , \quad t\in \left[0,1\right],
\end{equation*}
and, setting $s=\frac{1}{t} $ for $t\in \left(0,1\right]$, we have
\begin{equation*}
f^{\prime }\left(\frac{1}{s}\right)=ps^{1-p}\left[ \left(s+1\right)^{p-1}-\left(s-1\right)^{p-1}-2\right], \quad s\in [1,\infty ).
\end{equation*}
Put $g\left(s\right)=\left(s+1\right)^{p-1}-\left(s-1\right)^{p-1}-2$ for $s\in [1,\infty )$. We have
\begin{equation*}
g^{\prime }\left(s\right)=\left(p-1\right)\left[ \left(s+1\right)^{p-2}-\left(s-1\right)^{p-2}\right]>0, \quad s\in [1,\infty ),
\end{equation*}
so the function $g$ is increasing. Since $g\left(1\right)=2^{p-1}-2>0$, it follows that $g$ is positive and, consequently, $f^{\prime}$ is positive in $\left(0,1\right]$. Since $f\left(0\right)=0$, we conclude that $f\left(t\right) \geq 0$ for all $t\in \left[0,1\right]$.

We omit analogous argument for proving the remaining part.
 $\square $\\

Now we shall prove our main results of this section. Here and in what follows $\Delta$ denotes, as previously, the set 
$\{(r, \varepsilon)\in (0, \infty)^2: \varepsilon <2r \}$.

\begin{theorem}\label{thm1}
Let $\varphi :[0,\infty)\rightarrow [0,\infty )$ be an increasing bijection such that $\mathbf{p}_{\varphi }$ is a paranorm in 
$S=S(\Omega ,\Sigma ,\mu )$. If $\varphi $ is superquadratic, then the space 
$\left(S,\mathbf{p}_{\varphi }\right) $ is uniformly convex, and the function $\delta: \Delta\rightarrow (0, \infty)$, given by 
\begin{equation}\label{eA}
\delta \left(r,\varepsilon \right) =r-\varphi ^{-1}\left(\varphi \left(r\right) -\varphi \left(\frac{\varepsilon }{2}\right)\right),
\end{equation}
is its modulus of the convexity.
\end{theorem}

{\bf Proof.} 
Take an arbitrary $(r,\varepsilon) \in \Delta$ and then $x,y\in S $ such that 
\begin{equation*}
\mathbf{p}_{\varphi }\left(x\right) \leq r, \quad \mathbf{p}_{\varphi}\left(y\right) \leq r, \quad 
\text{and \ }\mathbf{p}_{\varphi }\left(x-y\right) \geq \varepsilon .
\end{equation*}
Putting 
\begin{equation*}
x_{1}=\frac{x+y}{2}, \qquad y_{1}=\frac{x-y}{2}
\end{equation*}
we have
\begin{equation*}
x=x_{1}+y_{1}, \qquad y=x_{1}-y_{1},
\end{equation*}
and, as $\varphi $ is increasing,
\begin{equation*}
\varphi \left(\mathbf{p}_{\varphi }\left(x_{1}+y_{1}\right) \right)+\varphi \left(\mathbf{p}_{\varphi }\left(x_{1}-y_{1}\right) \right) 
\leq 2\varphi \left(r\right).
\end{equation*}
Moreover, from the subadditivity of $\mathbf{p}_{\varphi }$ we have 
\begin{equation*}
\varepsilon \leq \mathbf{p}_{\varphi }\left(x-y\right) =\mathbf{p}_{\varphi
}\left(2y_{1}\right) \leq 2\mathbf{p}_{\varphi }\left(y_{1}\right) ,
\end{equation*}
that is 
\begin{equation*}
\text{\ }\mathbf{p}_{\varphi }\left(y_{1}\right) \geq \frac{\varepsilon }{2}.
\end{equation*}
From Lemma \ref{lem4} we get
\begin{equation*}
2\varphi \left(\mathbf{p}_{\varphi }\left(x_{1}\right) \right) +2\varphi\left(\mathbf{p}_{\varphi }\left(y_{1}\right) \right) 
\leq \varphi \left(\mathbf{p}_{\varphi }\left(x_{1}+y_{1}\right) \right) 
+\varphi \left(\mathbf{p}_{\varphi }\left(x_{1}-y_{1}\right) \right) \leq 2\varphi \left(r\right),
\end{equation*}
or
\begin{equation*}
\varphi \left(\mathbf{p}_{\varphi }\left(x_{1}\right) \right) +\varphi\left(\mathbf{p}_{\varphi }\left(y_{1}\right) \right) 
\leq \varphi \left(r\right).
\end{equation*}
Hence, as $\varphi $ is increasing, we have 
\begin{equation*}
\varphi \left(\mathbf{p}_{\varphi }\left(x_{1}\right) \right) \leq \varphi\left(r\right) 
-\varphi \left(\mathbf{p}_{\varphi }\left(y_{1}\right)\right) \leq \varphi \left(r\right) -\varphi \left(\frac{\varepsilon }{2}\right),
\end{equation*}
and thus
\begin{equation*}
\mathbf{p}_{\varphi }\left(\frac{x+y}{2}\right) 
=\mathbf{p}_{\varphi}\left(x_{1}\right) \leq \varphi ^{-1}\left(\varphi \left(r\right)-\varphi \left(\frac{\varepsilon }{2}\right) \right) 
=r-\delta(r,\varepsilon),
\end{equation*}
which completes the proof.
 $\square $

\medskip

From this result and Lemma \ref{lem1} we obtain

\begin{corollary}\label{thm2}
Assume that $\mu (\Omega )\leq 1. $ If $\varphi :[0,\infty )\rightarrow [0,\infty )$ is 
an increasing bijection and the function $F:[0, \infty)^2\rightarrow [0, \infty)$, given by \eqref{e1}, is concave, then $\mathbf{p}_{\varphi }$ is a paranorm in $S=S(\Omega,\Sigma ,\mu ) $. If, in addition, $\varphi $ is superquadratic, then the space $\left(S,\mathbf{p}_{\varphi }\right) $ is uniformly convex, and the function $\delta: \Delta \rightarrow (0, \infty)$, given by \eqref{eA}, is its modulus of convexity.
\end{corollary}

This result and Remark \ref{rem3} imply the following

\begin{corollary}\label{cor1}
Assume that $\mu (\Omega )\leq 1$. Let $\varphi :[0,\infty )\rightarrow [0,\infty )$ be a 
twice differentiable function with $\varphi(0)=0$, $\varphi ^{\prime }(r)>0 $ and $\varphi ^{\prime \prime}(r)>0 $ for all $r>0$ and such that the function $\frac{\varphi ^{\prime }}{\varphi ^{\prime\prime }}$ is superadditive in $(0, \infty)$. Then $\left(S,\mathbf{p}_{\varphi }\right) $ is a paranormed space. If, in addition, $\varphi $ is superquadratic, then the space $\left(S,\mathbf{p}_{\varphi }\right) $ is uniformly convex, and the function $\delta: \Delta \rightarrow (0, \infty)$, given by \eqref{eA}, is its modulus of convexity.
\end{corollary}

\begin{example}{\rm
For a fixed $p>2 $ consider the function $\varphi_p :[0,\infty )\rightarrow [0,\infty )$ given by 
$\varphi_p \left(t\right)=\frac{t^{p}}{t+1}.$ We have 
\begin{eqnarray*}
\frac{\varphi_p^{\prime }\left(t\right) }{\varphi_p^{\prime \prime }\left(t\right) }
=t\frac{\left(t+1\right)\left[\left(p-1\right)t+p\right]}{\left(p-1\right)\left(p-2\right)t^{2}+2p\left(p-2\right)t+p\left(p-1\right)},\quad t\geq 0.
\end{eqnarray*}
A standard, although tedious, computation shows that the function $t\mapsto \frac{\varphi_p^{\prime }\left(t\right) }{\varphi_p^{\prime \prime}\left(t\right)}/t$ is increasing 
in $\left(0,\infty \right)$. Hence we infer that  the function $\frac{\varphi_p^{\prime }}{\varphi_p^{\prime \prime }}$ is superadditive in $(0, \infty)$. Moreover, if $p=3 $ then for all\ $r\geq s\geq 0$ we obtain
\begin{eqnarray*}
&& \hspace{-1cm}\varphi_3 \left(r+s\right) +\varphi_3 \left(r-s\right) -2\varphi_3 \left(r\right) -2\varphi_3 \left(s\right) \\
&& \hspace{2cm}
=\frac{2s^{2}\left[ \left(r^{2}-s^{2}\right) \left(r+1\right)+r-s+2r\left(r+1\right) \right] }
				{\left(r+1\right)\left(s+1\right)\left(r+s+1\right) \left(r-s+1\right) }\geq 0,
\end{eqnarray*}
which shows that $\varphi_3$ is superquadratic. Thus there are non-power functions satisfying the assumptions of Corollary \ref{cor1}.
}\end{example}

If the measure space $(\Omega ,\Sigma ,\mu ) $ satisfies condition (ii) formulated in the Introduction, then, applying Lemma \ref{lem3} and Theorem \ref{thm1}, we obtain the following

\begin{corollary}\label{thm3}
Assume that for every $A\in \Sigma$ either $\mu (A)=0$, or $\mu (A)\geq 1$. Let $\varphi:[0,\infty )\rightarrow [0,\infty )$ be an increasing bijection, convex, and geometrically convex. Then $\mathbf{p}_{\varphi }$ is a paranorm in $S=S(\Omega ,\Sigma ,\mu )$. If, in addition, $\varphi $ is superquadratic, then the space 
$\left(S,\mathbf{p}_{\varphi }\right)$ is uniformly convex, and the function  $\delta: \Delta \rightarrow (0, \infty)$, given by \eqref{eA}, is its modulus of convexity.
\end{corollary}

\begin{example}{\rm
The function $\varphi:[0,\infty )\rightarrow [0,\infty )$, given by $\varphi \left(t\right) =t^{2}{\rm e}^{t}$, is an increasing bijection, convex and geometrically convex (cf. Remark \ref{rem5}). Setting 
\begin{equation*}
f\left(r,s\right)=\left(r+s\right) ^{2}{\rm e}^{r+s}+\left(r-s\right)^{2}{\rm e}^{r-s}-2r^{2}{\rm e}^{r}-2s^{2}{\rm e}^{s}, \quad 
r\geq s\geq 0,
\end{equation*}
we have $f\left(0,0\right) =0$. Making some calculations one can check that the derivatives $\partial_1f$ and $\partial_2f$ are nonnegative, so the function $f$ is increasing with respect to each variable. It follows that $f\left(r,s\right) \geq0 $ for all 
$r\geq s\geq 0$, that is the function $\varphi $ is superquadratic. Consequently, $\varphi$ satisfies the assumption of Corollary \ref{thm3}.
}\end{example}

To ensure the uniform convexity of the space $\left(S(\Omega ,\Sigma ,\mu ),\mathbf{p}_{\varphi }\right)$, up to now we have assumed that the function $\varphi$ is superquadratic. Now we show that this condition  may be replaced by another one, imposed on the function $H$ considered in Remark \ref{rem4}. We start with the following simple observation. 

\begin{lemma}\label{lem6}
For every function $\varphi: [0, \infty)\rightarrow \mathbb R$ the equality
$$
\varphi\left(\left|s\right|+\left|t\right|\right)+\varphi\left(\left|\left|s\right|-\left|t\right|\right|\right)=\varphi\left(\left|s+t\right|\right)+\varphi\left(\left|s-t\right|\right)
$$
holds for all $s,t \in \mathbb R$.
\end{lemma}

{\bf Proof.} It is enough to notice that if $st \geq 0$ then 
$$
\left|s\right|+\left|t\right|= \left|s+t\right|, \qquad   \left|\left|s\right|-\left|t\right|\right|=\left|s-t\right|
$$
and if $st\leq 0$ then
$$
\left|s\right|+\left|t\right|= \left|s-t\right|, \qquad   \left|\left|s\right|-\left|t\right|\right|=\left|s+t\right|
$$
for all $s,t \in \mathbb R$.  $ \square$
 
\begin{lemma}\label{lem5}
Let $\varphi :[0,\infty )\rightarrow [0,\infty )$ be a bijection with $\varphi \left(0\right) =0$ and 
$H:[0,\infty )^{2}\rightarrow [0,\infty )$ be given by 
\begin{equation}\label{eE}
H\left(r,s\right)=\varphi \left(\varphi ^{-1}\left(r\right) 
+\varphi^{-1}\left(s\right) \right) +\varphi \left(\left\vert \varphi ^{-1}\left(r\right) -\varphi ^{-1}\left(s\right)\right\vert \right).
\end{equation}
Then
\begin{eqnarray}\label{eB}
&& \hspace{0cm}\varphi \left(\mathbf{p}_{\varphi }\left(x\right) +\mathbf{p}_{\varphi}\left(y\right) \right) 
+\varphi \left(\left\vert \mathbf{p}_{\varphi}\left(x\right) -\mathbf{p}_{\varphi }\left(y\right) \right\vert \right)\\
&& \hspace{6cm} \leq \varphi \left(\mathbf{p}_{\varphi }\left(x+y\right) \right) +\varphi\left(\mathbf{p}_{\varphi }\left(x-y\right) \right) \nonumber
\end{eqnarray}
for all $x,y \in S(\Omega ,\Sigma ,\mu )$ iff the function $H$ satisfies the condition 
\begin{equation}\label{eC}
H\left(\sum^{k}_{i=1}a_i\left(r_i,s_i\right)\right)\leq \sum^{k}_{i=1}a_iH\left(r_i,s_i\right)
\end{equation}
for all $k \in {\mathbb N}, r_1, \ldots, r_k,s_1, \ldots s_k \in [0, \infty)$, and $a_1, \ldots, a_k \in [0, \infty)$ such that $a_i=\mu\left(A_i\right)$, $i =1, \ldots, k$, for some pairwise disjoint sets $A_1, \ldots, A_k \in \Sigma$.
\end{lemma}

{\bf Proof.} 
Assume that inequality \eqref{eB} holds for all $x,y \in S(\Omega ,\Sigma ,\mu )$. Take any numbers $r_1, \ldots, r_k$, $s_1, \ldots s_k \in [0, \infty)$, and  pairwise disjoint sets  $A_1, \ldots$, $A_k  \in \Sigma$ of finite measure; put $a_i=\mu \left(A_i\right)$, $i =1, \ldots, k$. Making use of \eqref{eB} for the functions 
\begin{equation*}
x=\sum_{i=1}^{k}\varphi^{-1}\left(r_{i}\right)\chi _{A_{i}}, \qquad y=\sum_{i=1}^{k}\varphi^{-1}\left(s_{i}\right)\chi _{A_{i}}
\end{equation*}
we obtain 
\begin{eqnarray*}
&& \hspace{-1cm} \varphi \left(\varphi ^{-1}\left(\sum_{i=1}^{k} a_{i}r_{i}\right) 
+\varphi ^{-1}\left(\sum_{i=1}^{k}a_{i}s_{i}\right) \right) \\
&& \hspace{4,5cm} +\varphi \left(\left\vert \varphi ^{-1}\left(\sum_{i=1}^{k}a_{i}r_{i}\right) 
-\varphi ^{-1}\left(\sum_{i=1}^{k}a_{i}s_{i} \right)\right\vert\right)
\\
&& \hspace{1cm}  \leq \sum_{i=1}^{k}a_i\left[ \varphi \left(\varphi ^{-1}\left(r_i\right)+ \varphi ^{-1}\left(s_i\right)\right) 
+\varphi\left(\left|\varphi ^{-1}\left(r_i\right)- \varphi ^{-1}\left(s_i\right) \right|\right) \right],
\end{eqnarray*}
i.e.
\begin{eqnarray*}
H\left(\sum_{i=1}^{k} a_{i}r_{i},\sum_{i=1}^{k} a_{i}s_{i}\right)\leq \sum_{i=1}^{k} a_{i}H\left(r_{i},s_i\right),
\end{eqnarray*}
which is \eqref{eC}. 

To prove the converse assume \eqref{eC} for all $k \in \mathbb N$, $r_1, \ldots, r_k$, $s_1, \ldots s_k \in [0, \infty)$, and $a_1, \ldots a_k \in [0, \infty)$ being the measures of pairwise disjoint sets. In other words we have the condition  
\begin{eqnarray}\label{eD}
H\left(\int_{\Omega}x{\rm d}\mu, \int_{\Omega}y{\rm d}\mu\right)\leq \int_{\Omega}H\circ\left(x,y\right){\rm d}\mu, 
\quad x,y \in S_{+}(\Omega ,\Sigma ,\mu ). 
\end{eqnarray}
Now take arbitrary $x,y \in S(\Omega ,\Sigma ,\mu )$. Then, by \eqref{eD} and Lemma \ref{lem6}, we get 
\begin{eqnarray*}
&& \hspace{-0,90cm}\varphi \left(\mathbf{p}_{\varphi }\left(x\right) +\mathbf{p}_{\varphi}\left(y\right) \right) 
+\varphi \left(\left\vert \mathbf{p}_{\varphi}\left(x\right) -\mathbf{p}_{\varphi }\left(y\right) \right\vert \right)\\
&& \hspace{1,25cm} =H\left(\int_{\Omega}\varphi\circ|x|{\rm d}\mu, \int_{\Omega}\varphi\circ|y|{\rm d}\mu\right)\leq \int_{\Omega}H\circ\left(\varphi\circ|x|,\varphi\circ|y|\right){\rm d}\mu\\
&& \hspace{1,25cm} =\int_{\Omega}\left[\varphi\left(\left|x(t)\right|+\left|y(t)\right|\right)+
\varphi\left(\left|\left|x(t)\right|-\left|y(t)\right|\right|\right)\right]{\rm d}\mu \\
&& \hspace{1,25cm} =\int_{\Omega}\left[\varphi\left(\left|x(t)+y(t)\right|\right)+\varphi\left(\left|x(t)-y(t)\right|\right)\right]
{\rm d}\mu \\
&& \hspace{1,25cm} =\int_{\Omega}\varphi\circ\left(\left|x+y\right|\right){\rm d}\mu 	      
                            +\int_{\Omega}\varphi\circ\left(\left|x-y\right|\right){\rm d}\mu \\
&& \hspace{1,25cm}=\varphi \left(\mathbf{p}_{\varphi }\left(x+y\right) \right) +\varphi\left(\mathbf{p}_{\varphi }\left(x-y\right) \right),
\end{eqnarray*}
that is \eqref{eB}.  $\square $

\begin{remark}\label{rem111}
Let $\varphi:[0, \infty)\rightarrow [0, \infty)$ be a bijection with $\varphi(0)=0$ and $H:[0, \infty)^2\rightarrow [0, \infty)$ be given by \eqref{eE}. If either \\

{\rm (i)} $\mu\left(\Omega\right)\leq 1$ and $H$ is convex,\\
or\\

{\rm (ii)} $\mu$ is the counting measure and $H$ is subadditive, \\
then condition \eqref{eC} holds for all $k \in \mathbb N$, $r_1, \ldots, r_k$, 
$s_1, \ldots, s_k\in [0, \infty)$, and $a_1, \ldots, a_k\in [0, \infty)$ such that $a_i=\mu\left(A_i\right)$, $i=1, \ldots, k$, for some pairwise disjoint sets $A_1, \ldots, A_k \in \Sigma$.
\end{remark}

{\bf Proof.} Assume (i) and take arbitrary  $r_1, \ldots, r_k$, $s_1, \ldots, s_k\in [0, \infty)$, and $a_1, \ldots, a_k\in [0, \infty)$ with $a_i=\mu\left(A_i\right)$, $i=1, \ldots, k$, where $A_1, \ldots, A_k \in \Sigma$ are pairwise disjoint. Then 
\begin{eqnarray*}
\sum^{k}_{i=1}a_i=\sum^{k}_{i=1}\mu \left(A_i\right)=\mu \left(\bigcup^{k}_{i=1}A_i\right)\leq 1.
\end{eqnarray*} 
Putting $r_{k+1}=s_{k+1}=0$ and $a_{k+1}=1-\sum^{k}_{i=1}a_i$, and making use of the convexity of $H$ we obtain 
\begin{eqnarray*}
H\left(\sum^{k}_{i=1}a_i\left(r_i,s_i\right)\right)=H\left(\sum^{k+1}_{i=1}a_i\left(r_i,s_i\right)\right)\leq \sum^{k+1}_{i=1}a_iH\left(r_i,s_i\right)=\sum^{k}_{i=1}a_iH\left(r_i,s_i\right),
\end{eqnarray*} 
as $H(0,0)=0$, and thus we come to \eqref{eC}.

If (ii) is satisfied, then the subadditivity of $H$ and a simple induction yields \eqref{eC} for all $r_1, \ldots, r_k$, $s_1, \ldots, s_k\in [0, \infty)$, and $a_1, \ldots, a_k\in {\mathbb N}_0$; moreover, since $\mu$ is the counting measure, we have $\mu\left(\Sigma\right) \subset {\mathbb N}_0 \cup \{\infty\}$, and the assertion follows. $\square$

\medskip 

Now we can prove

\begin{theorem}\label{thm4}
Let $\varphi :[0,\infty )\rightarrow [0,\infty )$ be an increasing bijection  such that $\mathbf{p}_{\varphi}$ is a paranorm in $S=S\left(\Omega, \Sigma, \mu\right)$. Assume that $\varphi$ is strictly convex and the function 
$H:[0,\infty)^{2}\rightarrow (0,\infty )$, given by \eqref{eE}, satisfies condition \eqref{eC} for all $k \in \mathbb N$, $r_1, \ldots, r_k$, $s_1, \ldots, s_k\in [0, \infty)$, and $a_1, \ldots, a_k\in [0, \infty)$ such that $a_i=\mu\left(A_i\right)$, $i=1, \ldots, k$, for some pairwise disjoint sets $A_1, \ldots, A_k \in \Sigma$. Then the space $\left(S, \mathbf{p}_{\varphi} \right)$ is uniformly convex and there is a modulus $\delta: \Delta \rightarrow [0,\infty )$ of convexity of it such that 
\begin{equation}\label{eF}
\varphi \left(r-\delta(r, \varepsilon) +\frac{\varepsilon }{2}\right) 
+\varphi \left(\left\vert r-\delta(r, \varepsilon) -\frac{\varepsilon }{2}\right\vert \right) =2\varphi\left(r\right) 
\end{equation}
for all $(r, \varepsilon) \in \Delta$.
\end{theorem}

{\bf Proof.} 
Take arbitrary $r \in(0, \infty)$, $\varepsilon \in \left(0,2r\right)$ and $x,y\in S$ such that 
\begin{equation*}
\mathbf{p}_{\varphi }\left(x\right) \leq r, \quad \mathbf{p}_{\varphi}\left(y\right) \leq r \quad {\rm and \ } \quad 
\mathbf{p}_{\varphi }\left(x-y\right) \geq \varepsilon .
\end{equation*}
Putting 
\begin{equation*}
x_{1}=\frac{x+y}{2},\qquad y_{1}=\frac{x-y}{2},
\end{equation*}
we have 
\begin{equation*}
x=x_{1}+y_{1},\qquad y=x_{1}-y_{1}.
\end{equation*}

The function $\lambda :\left[ 0,\infty \right)^{2}\rightarrow \left[ 0,\infty \right)$, defined by  
\begin{equation*}
\lambda \left(u,v\right) =\varphi \left(u+v\right) +\varphi \left(\left\vert u-v\right\vert \right),
\end{equation*}
is continuous. We prove that it is strictly increasing with respect to each of the variables. Since $\lambda$ is symmetric, it is enough to show that it is strictly increasing in first variable. Let us fix arbitrarily $v>0$. Clearly, the function $\lambda\left(\cdot,v\right)|_{[v, \infty)}$ is strictly increasing. To show that so is $\lambda\left(\cdot,v\right)|_{[0,v)}$ 
take any $u_1, u_2 \in \left[ 0,v\right] $, $u_1< u_2$. By the strict convexity of $\varphi $ we have
\begin{equation*}
\frac{\varphi \left(v-u_1\right)-\varphi \left(v-u_2\right)}{\left(v-u_1\right)-\left(v-u_2\right)}
<\frac{\varphi \left(v+u_2\right)-\varphi \left(v+u_1\right)}{\left(v+u_2\right)-\left(v+u_1\right)},
\end{equation*}
that is
\begin{equation*}
\varphi \left(v-u_1\right)-\varphi \left(v-u_2\right)<\varphi \left(v+u_2\right)-\varphi \left(v+u_1\right),
\end{equation*}
whence $\lambda\left(u_1,v\right)<\lambda\left(u_2,v\right)$.

Since $\frac{\varepsilon}{2}\leq \frac{1}{2}\mathbf{p}_{\varphi}(x-y) \leq \mathbf{p}_{\varphi}\left(\frac{x-y}{2}\right)
=\mathbf{p}_{\varphi}\left(y_1\right)$, it follows from the monotonicity of $\lambda\left(\mathbf{p}_{\varphi}\left(x_1\right),\cdot\right)$, Lemma \ref{lem5}, and the monotonicity of $\varphi$ that 
\begin{eqnarray*}
\lambda\left( \mathbf{p}_{\varphi }\left( x_{1}\right),\frac{\varepsilon}{2}\right)
&\leq & \lambda\left( \mathbf{p}_{\varphi }\left( x_{1}\right),\mathbf{p}_{\varphi}\left(y_{1}\right)\right) \\
&=& \varphi \left(\mathbf{p}_{\varphi }\left( x_{1}\right)+\mathbf{p}_{\varphi}\left( y_{1}\right)\right)
  +\varphi \left(\left\vert \mathbf{p}_{\varphi}\left(x_{1}\right)-\mathbf{p}_{\varphi }\left(y_{1}\right)\right\vert\right)\\ 
&\leq & \varphi \left(\mathbf{p}_{\varphi }\left( x_{1}+ y_{1}\right)\right)
  +\varphi \left(\mathbf{p}_{\varphi }\left( x_{1}-y_{1}\right)\right)\\
&=&\varphi \left(\mathbf{p}_{\varphi }\left(x\right)\right) +\varphi \left(\mathbf{p}_{\varphi }\left( y\right)\right)
\leq 2\varphi(r).
\end{eqnarray*}
Moreover,
\begin{eqnarray*}
\lambda \left(0,\frac{\varepsilon}{2}\right) =2\varphi\left(\frac{\varepsilon }{2}\right) 
<2\varphi \left(r\right)=\lambda(r,0)<\lambda \left(r,\frac{\varepsilon}{2}\right).
\end{eqnarray*} 
Hence, by the Darboux property of the continuous function $\lambda \left(\cdot,\frac{\varepsilon}{2}\right)$, there exists a unique 
$\delta\left(r,\varepsilon \right) \in \left(0,r\right)$ such that 
\begin{equation*}
\lambda \left(r-\delta\left(r, \varepsilon\right) ,\frac{\varepsilon }{2}\right) =2\varphi \left(r\right) .
\end{equation*}
The monotonicity of $\lambda $ implies that $\mathbf{p}_{\varphi }\left(x_{1}\right) <r-\delta\left(r, \varepsilon\right)$, that is 
\begin{equation*}
\mathbf{p}_{\varphi }\left(\frac{x+y}{2}\right) <r-\delta\left(r, \varepsilon\right).
\end{equation*}
Moreover, from the definition of the function $\lambda $ we get
\begin{equation*}
\varphi \left(r-\delta\left(r, \varepsilon\right) +\frac{\varepsilon }{2}\right)
+\varphi \left(\left\vert r-\delta\left(r, \varepsilon\right) -\frac{\varepsilon }{2}\right\vert \right) =2\varphi\left(r\right),
\end{equation*}
which completes the proof.
 $\square $

\medskip

Theorem \ref{thm4}, Lemma \ref{lem1}(b) and Remark \ref{rem111} imply the following
 
\begin{corollary}\label{cor2}
Assume that $\mu\left(\Omega\right)\leq 1$. If $\varphi :[0,\infty )\rightarrow [0,\infty )$ is an increasing bijection and the function $F:[0,\infty)^{2}\rightarrow [0,\infty )$, given by \eqref{e1}, is concave, then $\mathbf{p}_{\varphi}$ is a paranorm in $S=S\left(\Omega, \Sigma, \mu\right)$. If, in addition, $\varphi$ is strictly convex and the function $H:[0,\infty)^{2}\rightarrow (0,\infty )$, given by \eqref{eE}, is convex, then the space $\left(S, \mathbf{p}_{\varphi}\right)$ is uniformly convex and there is a modulus $\delta: \Delta\rightarrow (0, \infty)$ of convexity of it, satisfying equation \eqref{eF}.
\end{corollary}

On the other hand, making use of Theorem \ref{thm4}, Lemma \ref{lem3} and Remark \ref{rem111}, we come to 

\begin{corollary}\label{cor3}
Assume that $\mu $ is the counting measure on a set $\Omega \subset \mathbb N$. If $\varphi :[0,\infty )\rightarrow [0,\infty )$ is an increasing bijection, convex, and geometrically convex, then $\mathbf{p}_{\varphi}$ is a paranorm in $S=S\left(\Omega, 2^{\Omega}, \mu\right)$. 
If, in addition, $\varphi$ is strictly convex and the function $H:[0,\infty)^{2}\rightarrow (0,\infty )$, given by \eqref{eE}, is subadditive, then the space $\left(S, \mathbf{p}_{\varphi}\right)$ is uniformly convex and there is a modulus $\delta: \Delta\rightarrow (0, \infty)$ of convexity of it, satisfying equation \eqref{eF}.
\end{corollary}

The fact below follows from the observation made in Remark \ref{rem4}.  

\begin{remark}\label{rem7}
If $\varphi :[0,\infty )\rightarrow [0,\infty )$ is given by $\varphi \left(t\right) =t^{p}$ with a $p\in \left(1,2\right),$ then the function $H $ is convex.
\end{remark}

\begin{remark}\label{rem8}{\rm 
If a paranormed space $\left(S(\Omega,\Sigma,\mu),\mathbf{p}_{\varphi}\right)$ is uniformly convex, then its completion 
$\left(\mathcal{S}^{\varphi }(\Omega ,\Sigma ,\mu ),\mathbf{p}_{\varphi}\right)$  
is also a paranormed uniformly convex  space (with the same modulus of convexity). Thus, taking 
$\varphi :[0,\infty )\rightarrow [0,\infty )$, given by $\varphi \left(t\right) =t^{p}$ with a $p\in [2, \infty)$, and making use of 
Remark \ref{rem6} and Theorem \ref{thm1}, we come to the uniform convexity of the space $L^{p}\left(\Omega ,\Sigma ,\mu \right)$. 
If $p\in (1,2)$ then the same can be obtained by Remark \ref{rem7} and Corollary \ref{cor2} provided $\mu(\Omega)\leq 1$. 
If $\mu(\Omega)> 1$ we can use Theorem \ref{thm4} instead. However, to be perfectly honest then, it should be noticed that it is not so easy to verify condition \eqref{eC} in that case. Summarizing we obtain the celebrated Clarkson theorem on the uniform convexity of the $L^{p}$ spaces {\rm (see \cite{A.Clarc})} as a very special case of our results.
}\end{remark}

\section{The case of finite $\Omega$} 

Let $\Omega$ be finite, say $\Omega =\left\{1, \ldots, k\right\}$. Then $S\left(\Omega, \Sigma, \mu\right)={\mathbb R}^k$ no matter what a measure is $\mu$. Now, for a bijection $\varphi :[0,\infty )\rightarrow [0,\infty )$ vanishing at 0, the functional $\mathbf{p}_{\varphi}$ takes the form
\begin{equation*}
\mathbf{p}_{\varphi}(x)=\varphi^{-1}\left(\sum^{k}_{i=1}a_i\varphi\left(\left|x_i\right|\right)\right),
\end{equation*}
where $a_i=\mu\left(\left\{i\right\}\right)$ for $i=1, \ldots, k$. It turns out that under suitably weak assumptions imposed on $\varphi$ we can prove the uniform convexity of the space $\left(S,\mathbf{p}_{\varphi}\right)$.

\begin{theorem}\label{thm11}
Let $\varphi :[0,\infty )\rightarrow [0,\infty )$ be a strictly convex bijection such that $\mathbf{p}_{\varphi}$ is a paranorm in ${\mathbb R}^k$. Then the space $\left({\mathbb R}^k, \mathbf{p}_{\varphi}\right)$ is uniformly convex.
\end{theorem}

{\bf Proof. } 
Suppose to the contrary that the assertion fails to be true. Then there would exist $r \in (0, \infty)$, $\varepsilon \in (0, 2r)$, and sequences $\left(x_n\right)_{n \in \mathbb N}$ and $\left(y_n\right)_{n \in \mathbb N}$ of points of ${\mathbb R}^k$ such that 
\begin{equation*}
\mathbf{p}_{\varphi}\left(x_n\right)\leq r, \,\, \mathbf{p}_{\varphi}\left(y_n\right)\leq r,\,\, \mathbf{p}_{\varphi}\left(x_n-y_n\right)\geq \varepsilon \,\,\,\, {\rm and} \,\,\,\, \mathbf{p}_{\varphi}\left(\frac{x_n+y_n}{2}\right)> r-\frac{r}{n}
\end{equation*}
for each $n \in \mathbb N$. Since the sequences $\left(x_n\right)_{n \in \mathbb N}$ and $\left(y_n\right)_{n \in \mathbb N}$ are bounded, then, passing to subsequences if necessary, we can assume that they converge to $x \in {\mathbb R}^k$ and $y \in {\mathbb R}^k$, respectively. Hence, letting $n$ tend to $ \infty$, we get 
\begin{equation*}
\mathbf{p}_{\varphi}\left(x\right)\leq r, \,\, \mathbf{p}_{\varphi}\left(y\right)\leq r,\,\, \mathbf{p}_{\varphi}\left(x-y\right)\geq \varepsilon \,\,\,\, {\rm and} \,\,\,\, \mathbf{p}_{\varphi}\left(\frac{x+y}{2}\right)\geq r.
\end{equation*}
It follows from the assumptions that $\varphi$ is strictly increasing, and thus, by Remark \ref{rem11}, the function $\varphi\circ\left|\cdot\right|$ is strictly convex. As $\mathbf{p}_{\varphi}\left(x-y\right)\geq \varepsilon >0$ we know that $x\not=y$. Therefore
\begin{eqnarray*}
\varphi\left(r\right)
&\leq & \varphi\left(\mathbf{p}_{\varphi}\left(\frac{x+y}{2}\right)\right)=\sum^{k}_{i=1}a_i\varphi\left(\left|\frac{x_i+y_i}{2}\right|\right)\\
&< & \sum^{k}_{i=1}a_i \frac{\varphi\left(\left|x_i\right|\right)+\varphi\left(\left|y_i\right|\right)}{2}
=\frac{1}{2}\left(\varphi\left(\mathbf{p}_{\varphi}(x)\right)+\varphi\left(\mathbf{p}_{\varphi}(y)\right)\right)\leq \varphi(r),
\end{eqnarray*}
a contradiction. $\square$

\vspace{0,2cm}

Unfortunately, Theorem \ref{thm11} says nothing on a possible modulus of convexity $\delta\colon \Delta\rightarrow (0, \infty)$ of the space  $\left({\mathbb R}^k,\mathbf{p}_{\varphi}\right)$. The problem of the determining its effective form is important, however in some special cases can be nontrivial even on the real plane. We show this by considering an example of a paranormed space $\left({\mathbb R}^2,\mathbf{p}_{\varphi}\right)$, which is uniformly convex by virtue of Theorem \ref{thm11}; however, Theorems \ref{thm1} and \ref{thm4} are not applicable in that case, so to get any information on a modulus of convexity of that space we need to proceed in a different way.

\begin{example}\label{ex3}{\rm
Let $\mu$ be the counting measure on the set $\Omega=\{1,2\}$: $\mu(\{1\})=\mu(\{2\})=1$. Then the space 
$S\left(\Omega, 2^{\Omega}, \mu\right)$ is simply the plane ${\mathbb R}^2$. Define $\varphi: [0, \infty)\rightarrow [0, \infty)$ by $\varphi(t)={\rm e}^t-1$. It follows from Remark \ref{rem5} and Lemma \ref{lem5} that ${\mathbf{p}}_\varphi$, given by 
\begin{eqnarray}\label{eqn1}
{\mathbf{p}}_\varphi(x)=\varphi^{-1}\left(\varphi\left(\left|x_1\right|\right)+\varphi\left(\left|x_2\right|\right)\right),
\end{eqnarray}
is a paranorm in ${\mathbb R}^2$. Clearly, $\varphi$ is a strictly convex bijection, so the uniform convexity of the space $\left({\mathbb R}^2, {\mathbf{p}}_\varphi\right)$ follows from Theorem \ref{thm11}.

Setting $r=2\log 2$ and $s=\frac{1}{2}\log 2$ we  have
\begin{eqnarray*}
\varphi\left(r+s\right)+\varphi\left(\left|r-s\right|\right)&=& \varphi\left(\frac{5}{2}\log2\right)+\varphi\left(\frac{3}{2}\log2\right)\\
&=&\left(2^{\frac{5}{2}}-1\right)+\left(2^{\frac{3}{2}}-1\right)
<(6-1)+\left(2^{\frac{3}{2}}-1\right) \\
&=&2\left(2^{2}-1\right)+2\left(2^{\frac{1}{2}}-1\right)  
= 2\varphi\left(r\right)+2\varphi\left(s\right),
\end{eqnarray*}
which shows that $\varphi$ is not superquadratic, and thus Theorem \ref{thm1} is not applicable when determining a modulus of convexity of $\left({\mathbb R}^2, {\mathbf{p}}_\varphi\right)$.

To show that we cannot use also Theorem \ref{thm4} we observe that the function 
$H:\left[0, \infty\right)^2\rightarrow \left[0, \infty\right)$, given by
$$
H(r,s)=\varphi\left(\varphi^{-1}(r)+\varphi^{-1}(s)\right)+\varphi\left(\left|\varphi^{-1}(r)-\varphi^{-1}(s)\right|\right),
$$
does not satisfy condition \eqref{eC} for $k=2$ and some specific $r_1, r_2, s_1, s_2 \in [0, \infty)$ and $a_1=a_2=\mu(\{1\})=\mu(\{2\})$. Indeed, taking $r_1=r_2=s_1=s_2=1$ we have 
\begin{eqnarray*}
 && \hspace{-1,75cm}H\left(a_1\left(r_1,s_1\right)+a_2\left(r_2,s_2\right)\right)=H\left(2,2\right)\\
&=& \varphi\left(\varphi^{-1}(2)+\varphi^{-1}(2)\right)+\varphi\left(\left|\varphi^{-1}(2)-\varphi^{-1}(2)\right|\right)
=\varphi\left(2\varphi^{-1}\left(2\right)\right)\\ 
&=&\varphi\left(2\log 3\right)=9-1=8>6=2\cdot 3=2(4-1)
=2\varphi\left(2\log 2\right)\\
&=& 2\varphi\left(2\varphi^{-1}\left(1\right)\right) 
=\varphi\left(\varphi^{-1}(1)+\varphi^{-1}(1)\right)+\varphi\left(\left|\varphi^{-1}(1)-\varphi^{-1}(1)\right|\right)\\
&=& 2H(1,1)=a_1H\left(r_1,s_1\right)+a_2H\left(r_2,s_2\right).
\end{eqnarray*}
}\end {example}

Therefore to determine a modulus of convexity of the space $\left({\mathbb R}^2,{\mathbf{p}}_{\varphi}\right)$ we have to find another  independent argument. The proof of the last theorem of the paper provides a possible one.  Here and in what follows $\Delta=\{(r, \varepsilon) \in (0,\infty)^2\colon \, \varepsilon <2r\}$ and $\Delta_{\varphi}=\{(r, \varepsilon) \in (0,\infty)^2\colon \, \varphi(\varepsilon) \leq 2\varphi(r)\}$. Observe that if $(r, \varepsilon) \in \Delta$ then $\left(r, \frac{\varepsilon}{2}\right) \in \Delta_{\varphi}$ for any increasing $\varphi \colon \left[0, \infty\right)\rightarrow \left[0, \infty\right)$.

\begin{theorem}\label{thm5}
Let ${\mathbf{p}}_{\varphi}$ be the paranorm in ${\mathbb R}^2$ of form \eqref{eqn1}, where $\varphi: [0, \infty)$ $\rightarrow [0, \infty) $ is defined by $\varphi(t)={\rm e}^t-1 $. Then the function  $\delta_0\colon \Delta_{\varphi} \rightarrow \left(0, \infty\right)$, given by 
\begin{eqnarray*}
\delta_0\left(r,\varepsilon\right)
=r-\varphi^{-1}\left(\varphi\left(\frac{x_{r,\varepsilon}+r}{2}\right)
+\varphi\left(\frac{\varphi^{-1}\left(\varphi(r)-\varphi\left(x_{r,\varepsilon}\right)\right)}{2}\right)\right),
\end{eqnarray*}
where for every $\left(r,\varepsilon\right)\in \Delta_{\varphi}$ the number $x_{r,\varepsilon} \in \left[0,r\right]$ is the unique solution of the equation 
\begin{eqnarray}\label{eqn2}
\varphi(t)-\varphi(r-t)=\varphi(r)-\varphi(\varepsilon),
\end{eqnarray}
is strictly increasing with respect to second variable. Moreover, the function $\delta\colon \Delta \rightarrow (0, \infty),$ defined by $\delta(r,\varepsilon)=\delta_0\left(r,\frac{\varepsilon}{4}\right)$, is a modulus of convexity of the space $\left({\mathbb R}^2,\mathbf{p}_{\varphi}\right)$.
\end{theorem}

{\bf Proof. } The argument is divided into some parts. 

{\bf Claim A } {\it The function $\delta_0$ strictly increases in second variable and
\begin{eqnarray}\label{eqn3}
\mathbf{p}_{\varphi}\left(\frac{x+y}{2}\right)\leq r-\delta_0(r, \varepsilon)
\end{eqnarray} 
for all $(r, \varepsilon)\in \Delta_{\varphi}$ and $x,y\in [0, \infty)^2$ satisfying $\mathbf{p}_{\varphi}(x)=\mathbf{p}_{\varphi}(y)=r$ and $\mathbf{p}_{\varphi}(x-y)=\varepsilon$.}

To prove Claim A  fix a pair $\left(r,\varepsilon\right) \in \Delta_{\varphi}$ and any points $x, y \in [0, \infty)^2$ satisfying $\mathbf{p}_{\varphi}(x)=\mathbf{p}_{\varphi}(y)=r$ and $\mathbf{p}_{\varphi}(x-y)=\varepsilon$. Since $\mathbf{p}_{\varphi}(x-y)>0$, we have $x\not= y$. Assume, for instance, that $x_1\leq y_1$. Then $x_1< y_1$ and $x_2>y_2$, and  the equality $\mathbf{p}_{\varphi}\left(x-y\right)=\varepsilon $ means
 $$
\varphi\left(y_1-x_1\right)+\varphi\left(x_2-y_2\right)=\varphi(\varepsilon),
 $$
that is
 $$
{\rm e}^{y_1-x_1}+{\rm e}^{x_2-y_2}={\rm e}^{\varepsilon}+1.
 $$
Putting $u={\rm e}^{x_1} $, $v={\rm e}^{y_1} $, and taking into account that
 $$
u+{\rm e}^{x_2}={\rm e}^{x_1}+{\rm e}^{x_2}=\varphi\left(x_1\right)+\varphi\left(x_2\right)+2=\varphi\left(r\right)+2={\rm e}^r+1
 $$
and
 $$
v+{\rm e}^{y_2}={\rm e}^{y_1}+{\rm e}^{y_2}=\varphi\left(y_1\right)+\varphi\left(y_2\right)+2=\varphi\left(r\right)+2={\rm e}^r+1
 $$
we see that 
 $$
\frac{v}{u}+\frac{\rho-u}{\rho-v}=\alpha,
 $$
where $\rho ={\rm e}^r+1 $ and $\alpha={\rm e}^{\varepsilon}+1 $. Clearly $2<\alpha< 2\left(\rho-1\right)$ and
 $$
1\leq u<v \leq \rho -1.
 $$
In particular, $(u,v) $ is a point of the arc $\Gamma_{\alpha} $ given by
 $$
\Gamma_{\alpha}=\left\{(s,t)\in (0,\rho)^2: \frac{t}{s}+\frac{\rho -s}{\rho -t}=\alpha \right\};
 $$
consequently, $\Gamma_{\alpha}\cap \left[1,\rho-1\right]^2\not= \emptyset $. Since $\alpha>2 $ we see that $\Gamma_{\alpha} $ has no common points with the diagonal, and thus, as $u<v $, $\Gamma_{\alpha} $ lies strictly over the diagonal. Defining $h_{\alpha}:\left(0, \rho \right)^2\rightarrow{\mathbb R} $ by
 $$
h_{\alpha}(s,t)=\frac{t}{s}+\frac{\rho-s}{\rho-t}-\alpha
 $$
we observe that $\Gamma_{\alpha} $ coincides with the set of zeros of $h_{\alpha} $. For all $s,t \in (0, \rho) $ we have
 $$
\partial_1h_{\alpha}(s,t)=-\frac{t}{s^2}-\frac{1}{\rho-t}<0
 $$
and
 $$
\partial_2h_{\alpha}(s,t)=\frac{1}{s}+\frac{\rho-s}{\left(\rho-t\right)^2}>0,
 $$
whence rank $h'(s,t)=1$. Moreover,
 $$
- \frac{\partial_1h_{\alpha}(s,t)}{\partial_2h_{\alpha}(s,t)}>0, \quad (s,t)\in \Gamma_{\alpha},
 $$
so making use of Implicit Function Theorem we infer that $\Gamma_{\alpha} $ is the graph of a strictly increasing function of class $C^1 $. Letting $s $ tend to 0 in the equality defining $\Gamma_{\alpha} $ we see that $(0,0) $ is one of the endpoints of $\Gamma_{\alpha} $. Similarly, when $t $ tends to $\rho $ there, we deduce that the other one is $\left(\rho,\rho\right)$.

We have
\begin{eqnarray}\label{eqn4}
\mathbf{p}_{\varphi}\left(\frac{x+y}{2}\right)&=&\varphi^{-1}\left(\varphi\left(\frac{x_1+y_1}{2}\right)+\varphi\left(\frac{x_2+y_2}{2}\right)\right) \\ &=&\log\left(\sqrt{uv}+\sqrt{\left(\rho-u\right)\left(\rho-v\right)}-1\right), \nonumber
\end{eqnarray}
and thus to prove \eqref{eqn3} we need to maximalize the function $f:\left\{(s,t)\in \left[0, \rho\right]^2:\right.$ $\left. s\leq t\right\}\rightarrow {\mathbb R} $, given by
 $$
f(s,t)=\sqrt{st}+\sqrt{\left(\rho-s\right)\left(\rho-t\right)},
 $$
on the arc $\Gamma_{\alpha}\cap \left[1,\rho-1\right]^2 $. Of course $f|_{{\rm cl} \Gamma_{\alpha}} $ takes its minimum and maximum values.

If $0<s<t< \rho $ then
 $$
\partial_1f(s,t)=\frac{t}{2\sqrt{st}}-\frac{\rho-t}{2\sqrt{\left(\rho-s\right)\left(\rho-t\right)}}=\frac{1}{2}\left(\sqrt{\frac{t}{s}}-\sqrt{\frac{\rho-t}{\rho-s}}\right)>0
 $$
and
 $$
\partial_2f(s,t)=\frac{s}{2\sqrt{st}}-\frac{\rho-s}{2\sqrt{\left(\rho-s\right)\left(\rho-t\right)}}=\frac{1}{2}\left(\sqrt{\frac{s}{t}}-\sqrt{\frac{\rho-s}{\rho-t}}\right)<0,
 $$
and thus $f $ is strictly increasing with respect to the first variable and strictly decreasing as a function of the second one. In particular, if $0\leq s<t \leq \rho $ then
 $$
f(s,t)<f(t,t)=\sqrt{t^2}+\sqrt{\left(\rho-t\right)^2}=t+\left(\rho-t\right)=\rho;
 $$
moreover, $f(s,s)=\rho $ for every $s \in \left[0, \rho\right] $. Therefore $f|_{{\rm cl} \Gamma_{\alpha}} $ takes its maximum value $\rho $ at the points $(0,0) $ and $(\rho, \rho) $.

Let $(s,t)\in \Gamma_{\alpha}$ be a point where $f|_{\Gamma_{\alpha}} $ takes a local extremum. Then there is a $\lambda \in {\mathbb R} $ such that 
 $$
\partial_1f(s,t)=\lambda\partial_1h_{\alpha}(s,t)
 $$
and 
 $$
\partial_2f(s,t)=\lambda\partial_2h_{\alpha}(s,t),
 $$
that is we solve the system consisting of the equations
\begin{eqnarray*}
\frac{1}{2}\left(\sqrt{\frac{t}{s}}-\sqrt{\frac{\rho-t}{\rho-s}}\right)&=&-\lambda\left(\frac{t}{s^2}+\frac{1}{\rho-t}\right),\\
\frac{1}{2}\left(\sqrt{\frac{s}{t}}-\sqrt{\frac{\rho-s}{\rho-t}}\right)&=&\lambda\left(\frac{1}{s}+\frac{\rho-s}{\left(\rho-t\right)^2}\right)
\end{eqnarray*}
and 
\begin{eqnarray}\label{ej3}
\frac{t}{s}+\frac{\rho-s}{\rho-t}=\alpha.
\end{eqnarray}
If $\lambda =0 $ then we would have $\sqrt{\frac{t}{s}}-\sqrt{\frac{\rho-s}{\rho-t}}=0 $, whence $\frac{t}{s}=\frac{\rho-s}{\rho-t} $, that is $t\left(\rho-s\right)=s\left(\rho-t\right) $ and, consequently, $s=t $ which is impossible. Thus $\lambda\not=0 $ and
\begin{eqnarray*}
\frac{\sqrt{\frac{s}{t}}-\sqrt{\frac{\rho-s}{\rho-t}}}{\sqrt{\frac{t}{s}}-\sqrt{\frac{\rho-t}{\rho-s}}}=
-\frac{\frac{1}{s}+\frac{\rho-s}{\left(\rho-t\right)^2}}{\frac{t}{s^2}+\frac{1}{\rho-t}}
=- \frac{s\left(\rho-t\right)^2+s^2\left(\rho-s\right)}{t\left(\rho-s\right)^2+s^2\left(\rho-t\right)}.
\end{eqnarray*}
Due to the identity
\begin{eqnarray*}
\frac{a-b}{\frac{1}{a}-\frac{1}{b}}=-ab
\end{eqnarray*}
held for any $a,b \in {\mathbb R}\setminus \{0\} $, $a\not=b $, the previous equalities give
\begin{eqnarray}\label{ej4}
\sqrt{\frac{s}{t}}\sqrt{\frac{\rho-s}{\rho-t}}
=\frac{s\left(\rho-t\right)^2+s^2\left(\rho-s\right)}{t\left(\rho-t\right)^2+s^2\left(\rho-t\right)}.
\end{eqnarray}
Making use of {\it Mathematica 4.0} one can check that the system of equations \eqref{ej3} and \eqref{ej4} has a unique solution $(s,t)$, viz. $\left(\frac{2}{\alpha+2}\rho,\frac{\alpha}{\alpha+2}\rho \right)$, satisfying the condition $0<s<t< \rho $; clearly $f\left(\frac{2}{\alpha+2}\rho,\frac{\alpha}{\alpha+2}\rho \right)<\rho=f(0,0)=f(\rho,\rho) $. Consequently, the function $f|_{\Gamma_{\alpha}} $ has only one local extremum, its absolute minimum taken at the point $\left(\frac{2}{\alpha+2}\rho,\frac{\alpha}{\alpha+2}\rho \right)$, and thus the absolute maximum $c\left(\rho,\alpha\right)$ of the function $f|_{\Gamma_{\alpha} \cap \left[1,\rho-1\right]^2}$ is taken at at least one of the endpoints of the arc ${\Gamma_{\alpha} \cap \left[1,\rho-1\right]^2}$. One of them is $\left(s_{\alpha},{\rho-1}\right)$ with some $s_{\alpha} \in \left[1,\rho-1\right)$. Since $\left(s_{\alpha},{\rho-1}\right) \in \Gamma_{\alpha}$, also $\left(1,\rho-s_{\alpha}\right)\in \Gamma_{\alpha}$. So the last is another endpoint of the arc ${\Gamma_{\alpha} \cap \left[1,\rho-1\right]^2}$. Observe that 
 $$
f\left(s_{\alpha},{\rho-1}\right)=\sqrt{s_{\alpha}\left(\rho-1\right)}+\sqrt{\rho-s_{\alpha}}=f\left(1,\rho-s_{\alpha}\right),
 $$
and thus the maximum value $c\left(\rho,\alpha\right)$ is taken at both the endpoints of the arc.

Setting $x_{r, \varepsilon}= \log s_{\alpha}$ we have $x_{r, \varepsilon} \in [0, r)$. Moreover, as $\left(s_{\alpha},{\rho-1}\right) \in \Gamma_{\alpha}$, we obtain
\begin{eqnarray*}
\varphi\left(x_{r, \varepsilon}\right)-\varphi\left(r-x_{r, \varepsilon}\right)
&=&{\rm e}^{x_{r, \varepsilon}}-{\rm e}^{r-x_{r, \varepsilon}}=s_{\alpha}-\frac{\rho-1}{s_{\alpha}}\\
&=&\rho-\alpha={\rm e}^r-{\rm e}^{ \varepsilon}=\varphi\left(r\right)-\varphi\left(\varepsilon\right).
\end{eqnarray*}
Since the function $[0,r]\ni t\mapsto \varphi\left(t\right)-\varphi\left(r-t\right)$ is strictly increasing, it follows that $x_{r, \varepsilon}$ is the unique solution of equation \eqref{eqn2}. Observe also that, by virtue of \eqref{eqn4}, we have 
\begin{eqnarray*}
\mathbf{p}_{\varphi}\left(\frac{x+y}{2}\right)
&=&\log\left(\sqrt{uv}+\sqrt{\left(\rho-u\right)\left(\rho-v\right)}-1\right)=\log\left(f(u,v)-1\right)\\
& \leq & \log\left(f\left(s_{\alpha},\rho-1\right)-1\right)=\log \left(\sqrt{s_{\alpha}\left(\rho-1\right)}+\sqrt{\rho-s_{\alpha}}-1\right)\\
&=&\log \left({\rm e}^{\frac{x_{r,\varepsilon}+r}{2}}+\sqrt{{\rm e}^r-{\rm e}^{x_{r,\varepsilon}}+1}-1\right)\\
&=&\varphi^{-1}\left(\varphi\left({\frac{x_{r,\varepsilon}+r}{2}}\right)
+\varphi\left(\frac{\varphi^{-1}\left(\varphi(r)-\varphi\left(x_{r,\varepsilon}\right)\right)}{2}\right)\right)\\
&=&r-\delta_0(r,\varepsilon).
\end{eqnarray*}
If $s_{{\alpha}_1} \leq s_{{\alpha}_2}$ for some ${{\alpha}_1}, {{\alpha}_2} \in \left(2, 2\left(\rho-1\right)\right]$ satisfying ${{\alpha}_1} < {{\alpha}_2}$, then, since $\left(s_{{\alpha}_1},\rho -1\right)\in \Gamma_{{\alpha}_1}$ and $\left(s_{{\alpha}_2},\rho-1 \right)\in \Gamma_{{\alpha}_2}$, we would have
\begin{eqnarray*}
\frac{\rho-1}{s_{{\alpha}_1}}+\rho-s_{{\alpha}_1}={{\alpha}_1}<{{\alpha}_2}=\frac{\rho-1}{s_{{\alpha}_2}}+\rho-s_{{\alpha}_2} \leq \frac{\rho-1}{s_{{\alpha}_1}}+\rho-s_{{\alpha}_1}
\end{eqnarray*}
which is impossible. Thus the function $\left(2, 2\left(\rho-1\right)\right]\ni \alpha \mapsto s_{\alpha}$ strictly decreases. Therefore, since 
$$
\delta_0\left(r, \varepsilon\right)=\delta_0\left(r, \log\left(\alpha-1\right)\right)=r-\log\left(f\left(s_{\alpha},\rho-1\right)-1\right)
$$
for each $\left(r,\varepsilon\right) \in \Delta_{\varphi}$ and $f$ strictly increases as a function of first variable, the function $\delta_0$ strictly increases in second variable. 
This completes the proof of Claim A. 

{\bf Claim B } {\it If $\left(r, \varepsilon\right) \in \Delta$ then
\begin{eqnarray*}
\mathbf{p}_{\varphi}\left(\frac{x+y}{2}\right)\leq r-\delta_0\left(r,\frac{\varepsilon}{2}\right)
\end{eqnarray*}
for all $x, y \in {\mathbb R}^2$ satisfying $\mathbf{p}_{\varphi}(x)=\mathbf{p}_{\varphi}(y)=r$, $\mathbf{p}_{\varphi}(x-y)\geq \varepsilon$, and 
\begin{eqnarray}\label{eqn5}
\left|x_1\right|\leq y_1 \qquad {\rm and} \qquad \left|y_2\right|\leq x_2.
\end{eqnarray}}

To see this fix any $(r,\varepsilon) \in \Delta$. Then $\left(r,\frac{\varepsilon}{2}\right) \in \Delta_\varphi$. Take also arbitrary points $x, y \in {\mathbb R}^2$ such that $\mathbf{p}_{\varphi}(x)=\mathbf{p}_{\varphi}(y)=r$, $\mathbf{p}_{\varphi}(x-y)\geq \varepsilon$, and \eqref{eqn5} holds.

\vspace{0,2cm}

If $x,y \in [0, \infty)^2$ then, setting $\varepsilon'= \mathbf{p}_{\varphi}(x-y)$, we have $\varepsilon'\geq \varepsilon$ and, by \eqref{eqn5}, 
\begin{eqnarray*}
\varphi\left(\varepsilon'\right)=\varphi\left(\mathbf{p}_{\varphi}(x-y)\right)=\varphi\left(y_1-x_1\right)+\varphi\left(x_2-y_2\right)\leq 2\varphi(r),
\end{eqnarray*}
that is $\left(r,\varepsilon' \right) \in \Delta_\varphi$. Thus Claim A implies 
\begin{eqnarray*}
\mathbf{p}_{\varphi}\left(\frac{x+y}{2}\right)\leq r-\delta_0\left(r,\varepsilon'\right)\leq r-\delta_0\left(r,\varepsilon\right)\leq r-\delta_0\left(r,\frac{\varepsilon}{2}\right).
\end{eqnarray*}

Now assume that exactly one of the points $x, y$, say $x$, lies in $[0, \infty)^2$. Then, in view of \eqref{eqn5}, we have 
\begin{eqnarray*}
0\leq x_1 \leq y_1 \qquad {\rm and} \qquad 0<-y_2\leq x_2.
\end{eqnarray*}
Let $x'=x$ and $y'=(r,0)$. Then $x', y' \in [0, \infty)^2$. Then $\mathbf{p}_{\varphi}\left(x'\right)=\mathbf{p}_{\varphi}\left(y'\right)=r$ and, putting $\varepsilon'=\mathbf{p}_{\varphi}\left(x'-y'\right)$, we get
\begin{eqnarray*}
\varphi\left(\varepsilon'\right)=\mathbf{p}_{\varphi}\left(x'-y'\right)=\varphi\left(r-x_1\right)+\varphi\left(x_2-0\right)\leq 2\varphi(r),
\end{eqnarray*}
whence $\left(r, \varepsilon'\right)\in \Delta_\varphi$. Moreover,
\begin{eqnarray*}
\varphi\left(\varepsilon'\right)&=&\varphi\left(r-x_1\right)+\varphi\left(x_2\right) 
\geq \varphi\left(\frac{y_1-x_1}{2}\right)+\varphi\left(\frac{x_2-y_2}{2}\right)\\
&=& \varphi\left(\mathbf{p}_{\varphi}\left(\frac{x-y}{2}\right)\right)
\geq \varphi\left(\frac{1}{2}\mathbf{p}_{\varphi}\left(x-y\right)\right) 
\geq \varphi\left(\frac{\varepsilon}{2}\right),
\end{eqnarray*}
that is $\varepsilon' \geq \frac{\varepsilon}{2}$, and thus, by Claim A, 
\begin{eqnarray*}
\mathbf{p}_{\varphi}\left(\frac{x'+y'}{2}\right) \leq r-\delta_0\left(r,\varepsilon' \right)
\leq r-\delta_0\left(r,\frac{\varepsilon}{2} \right).
\end{eqnarray*}
Now, since 
\begin{eqnarray*}
\varphi\left(\mathbf{p}_{\varphi}\left(\frac{x+y}{2}\right)\right)&=& \varphi\left(\frac{x_1+y_1}{2}\right)+\varphi\left(\frac{x_2+y_2}{2}\right)\\
&\leq & \varphi\left(\frac{{x}_1'+r}{2}\right)+\varphi\left(\frac{{x}_2'+0}{2}\right)
=\varphi\left(\mathbf{p}_{\varphi}\left(\frac{x'+y'}{2}\right)\right),
\end{eqnarray*}
it follows that
\begin{eqnarray*}
\mathbf{p}_{\varphi}\left(\frac{x+y}{2}\right)\leq \mathbf{p}_{\varphi}\left(\frac{x'+y'}{2}\right)
\leq r-\delta_0\left(r, \frac{\varepsilon}{2}\right).
\end{eqnarray*}

Finally consider the case $x,y \not\in [0,\infty)^2$. Then, by virtue of \eqref{eqn5}, we get
\begin{eqnarray*}
0< -x_1 \leq y_1 \qquad {\rm and} \qquad 0<-y_2\leq x_2.
\end{eqnarray*}
Let $x'=(0,r)$ and $y'=(r,0)$. Then $x', y' \in [0, \infty)^2$,  $\mathbf{p}_{\varphi}\left(x'\right)=\mathbf{p}_{\varphi}\left(y'\right)=r$ and 
\begin{eqnarray*}
\varphi\left(\mathbf{p}_{\varphi}\left(x'-y'\right)\right)&=& \varphi(r)+\varphi(r)
\geq \varphi\left(\frac{y_1+y_1}{2}\right)+\varphi\left(\frac{x_2+x_2}{2}\right)\\
&\geq & \varphi\left(\frac{y_1-x_1}{2}\right)+\varphi\left(\frac{x_2-y_2}{2}\right)
=\varphi\left(\mathbf{p}_{\varphi}\left(\frac{x-y}{2}\right)\right)\\
&\geq & \varphi\left(\frac{1}{2}\mathbf{p}_{\varphi}\left(x-y\right)\right)
\geq \varphi\left(\frac{\varepsilon}{2}\right),
\end{eqnarray*}
i.e. $\mathbf{p}_{\varphi}\left(x'-y'\right) \geq \frac{\varepsilon}{2}$, whence again by Claim A, 
\begin{eqnarray*}
\mathbf{p}_{\varphi}\left(\frac{x'+y'}{2}\right) \leq r-\delta_0\left(r,\frac{\varepsilon}{2} \right);
\end{eqnarray*}
consequently, 
\begin{eqnarray*}
\varphi\left(\mathbf{p}_{\varphi}\left(\frac{x+y}{2}\right)\right)&=& \varphi\left(\frac{x_1+y_1}{2}\right)+\varphi\left(\frac{x_2+y_2}{2}\right)\\
&\leq & \varphi\left(\frac{0+r}{2}\right)+\varphi\left(\frac{r+0}{2}\right)
=\varphi\left(\mathbf{p}_{\varphi}\left(\frac{x'+y'}{2}\right)\right),
\end{eqnarray*}
whence
\begin{eqnarray*}
\mathbf{p}_{\varphi}\left(\frac{x+y}{2}\right) \leq \mathbf{p}_{\varphi}\left(\frac{x'+y'}{2}\right)
\leq r-\delta_0\left(r, \frac{\varepsilon}{2}\right).
\end{eqnarray*}
This completes the proof of Claim B.

To prove the theorem fix arbitrarily $\left(r,\varepsilon\right) \in \Delta$, and then points $x,y \in {\mathbb R}^2 $ satisfying the inequalities $\mathbf{p}_{\varphi}(x)\leq r $, $\mathbf{p}_{\varphi}(y)\leq r $ and $\mathbf{p}_{\varphi}(x-y)\geq \varepsilon $. Observe that the ball
 $$
\left\{(s,t)\in {\mathbb R}^2: \mathbf{p}_{\varphi}(s,t)\leq r\right\}=\left\{(s,t)\in {\mathbb R}^2: \varphi\left(\left|s\right|\right)+\varphi\left(\left|t\right|\right)\leq\varphi(r)\right\}
 $$ 
is symmetric with respect to the straight lines with the equations $t=0 $, $t=s $, $s=0 $ and $t=-s $. Therefore, without loss of generality, we may assume that the mean $\frac{x+y}{2} $ lies in the first quadrant of the real plane:
 $$
x_1+y_1\geq 0 \qquad {\rm and} \qquad x_2+y_2\geq 0.
 $$
Interchanging $x $ and $y $ we may additionally assume that $x_1\leq y_1 $. Thus 
\begin{eqnarray}\label{ej6}
\left|x_1\right|\leq y_1 \qquad {\rm and} \qquad x_2+y_2\geq 0.
\end{eqnarray}
Two complementary cases should be considered: (i) $y_2 \leq x_2 $, and (ii) $y_2> x_2 $.\\

(i) Then, by \eqref{ej6}, we have 
 $$
\left|x_1\right|\leq y_1 \qquad {\rm and} \qquad \left|y_2\right|\leq x_2.
 $$
Since $\mathbf{p}_{\varphi}(x)\leq r $, we have 
$\varphi\left(\left|x_1\right|\right)+\varphi\left(\left|x_2\right|\right)\leq \varphi(r) $, whence
\begin{eqnarray}\label{ej7}
x_2 \leq \pi_r\left(\left|x_1\right|\right),
\end{eqnarray}
where $\pi_r\colon [0,r]\rightarrow[0,r]$ is given by $\pi_r(t)=\varphi^{-1}\left(\varphi(r)-\varphi(t)\right)$.  
Similarly, $\mathbf{p}_{\varphi}\left(y\right)\leq r $ implies $\varphi\left(\left|y_1\right|\right)+\varphi\left(\left|y_2\right|\right)\leq \varphi(r) $ which gives 
\begin{eqnarray}\label{ej8}
y_1 \leq \pi_r\left(\left|y_2\right|\right).
\end{eqnarray}
Setting $x'=\left(x_1,\pi_r\left(\left|x_1\right|\right)\right) $ and $y'=\left(\pi_r\left(\left|y_2\right|\right),y_2\right)$ we see that $\mathbf{p}_{\varphi}\left(x'\right)=\mathbf{p}_{\varphi}\left(y'\right)=r$. It follows from \eqref{ej8} and \eqref{ej7} that 
 $$
x_1'=x_1\leq y_1 \leq y_1' \qquad {\rm and} \qquad x_2'\geq x_2\geq y_2=y_2'.
 $$
Hence
\begin{eqnarray*}
\varphi\left(\mathbf{p}_\varphi\left(x'-y'\right)\right)& =& \varphi\left(y_1'-x_1'\right)+\varphi\left(x_2'-y_2'\right)\\
&\geq & \varphi\left(y_1-x_1\right)+\varphi\left(x_2-y_2\right)
=\varphi\left(\mathbf{p}_\varphi\left(x-y\right)\right).
\end{eqnarray*}
Then, putting $\varepsilon'=\mathbf{p}_\varphi\left(x'-y'\right)$, we see that $\varepsilon'\geq \mathbf{p}_\varphi\left(x-y\right)\geq \varepsilon$. Moreover, 
 $$
x_1'+y_1'\geq x_1+y_1\geq 0, \quad x_2'+y_2'\geq x_2+y_2\geq 0
 $$ 
and 
\begin{eqnarray*}
\varphi\left(\mathbf{p}_\varphi\left(\frac{x+y}{2}\right)\right)&=&
\varphi\left(\frac{x_1+y_1}{2}\right)+\varphi\left(\frac{x_2+y_2}{2}\right) \\
&\leq & \varphi\left(\frac{x_1'+y_1'}{2}\right)+\varphi\left(\frac{x_2'+y_2'}{2}\right) 
=\varphi\left(\mathbf{p}_\varphi\left(\frac{x'+y'}{2}\right)\right),
\end{eqnarray*}
whence, applying Claim B to the points $x'$ and $y' $, we get 
\begin{eqnarray*}
\mathbf{p}_\varphi\left(\frac{x+y}{2}\right) &\leq& \mathbf{p}_\varphi\left(\frac{x'+y'}{2}\right)\leq r-\delta_0\left(r, \frac{\varepsilon'}{2}\right)\leq r-\delta_0\left(r, \frac{\varepsilon}{4}\right)\\
&=& r-\delta\left(r, \varepsilon\right).
\end{eqnarray*}

(ii) Now, by \eqref{ej6}, we come to
 $$
\left|x_1\right|\leq y_1 \qquad {\rm and} \qquad \left|x_2\right|\leq y_2.
 $$
Then, since the function $\pi_r $ decreases, we have 
\begin{eqnarray}\label{ej11}
y_1\leq \pi_r\left(y_2\right)\leq \pi_r\left(\left|x_2\right|\right)  \quad {\rm and} \quad y_2\leq \pi_r\left(y_1\right)\leq \pi_r\left(\left|x_1\right|\right).
\end{eqnarray}
Putting $x'=\left(x_1,\pi_r\left(\left|x_1\right|\right) \right) $, $x''=\left(\pi_r\left(\left|x_2\right|\right),x_2 \right) $, $y'=\left(y_1,\pi_r\left(y_1\right) \right) $ and $y''=\left(\pi_r\left(y_2\right),y_2 \right)$ we see that $\mathbf{p}_\varphi\left(x'\right)=\mathbf{p}_\varphi\left(x''\right)=\mathbf{p}_\varphi\left(y'\right)=\mathbf{p}_\varphi\left(y''\right)=r$. By \eqref{ej11} we have 
 $$
x_1'=x_1\leq y_1 \leq x_1'' \qquad {\rm and} \qquad x_2''=x_2 \leq y_2 \leq x_2',
 $$
whereas \eqref{ej8} and the inequality $y_2\leq\pi_r\left(y_1\right)$ give 
 $$
y_1=y_1'\leq y_1'' \qquad {\rm and} \qquad y_2=y_2'' \leq y_2',
 $$
respectively. 

Thus
\begin{eqnarray*}
\varphi\left(\varepsilon\right)
&\leq& \varphi\left(\mathbf{p}_\varphi\left(x-y\right)\right)
= \varphi\left(y_1-x_1\right)+\varphi\left(y_2-x_2\right)\\
&\leq& \varphi\left(x_1''-x_1'\right)+\varphi\left(x_2'-x_2''\right)
=\varphi\left(\mathbf{p}_\varphi\left(x'-x''\right)\right),
\end{eqnarray*}
whence
 $$
\varepsilon \leq \mathbf{p}_{\varphi}\left(x'-x''\right)\leq \mathbf{p}_{\varphi}\left(x'-\bar{y}\right)+\mathbf{p}_{\varphi}\left(x''-\bar{y}\right),
 $$
where $\bar{y} $ is chosen as either $y' $, or $y'' $. Assume, for instance, that $\mathbf{p}_{\varphi}\left(x''-\bar{y}\right)\geq \frac{\varepsilon}{2} $. Then, putting $\varepsilon'=\mathbf{p}_{\varphi}\left(x''-\bar{y}\right) $, we have $\varepsilon'\geq \frac{\varepsilon}{2} $. Moreover, 
 $$
x_1''+\bar{y_1}\geq x_1+y_1\geq 0, \quad x_2''+\bar{y_2}\geq x_2+y_2\geq 0
 $$ 
and
\begin{eqnarray*}
\varphi\left(\mathbf{p}_\varphi\left(\frac{x+y}{2}\right)\right)
&=& \varphi\left(\frac{x_1+y_1}{2}\right)+\varphi\left(\frac{x_2+y_2}{2}\right)\\
&\leq & \varphi\left(\frac{x_1''+\bar{y_1}}{2}\right)+\varphi\left(\frac{x_2''+\bar{y_2}}{2}\right) 
=\varphi\left(\mathbf{p}_\varphi\left(\frac{x''+\bar{y}}{2}\right)\right),
\end{eqnarray*}
whence, by Claim B, we obtain
 $$
\mathbf{p}_\varphi\left(\frac{x+y}{2}\right)\leq \mathbf{p}_\varphi\left(\frac{x''+\bar{y}}{2}\right) 
\leq r-\delta_0\left(r,\frac{\varepsilon'}{2}\right)\leq r-\delta_0\left(r,\frac{\varepsilon}{4}\right)
\leq r-\delta\left(r,\varepsilon\right),
 $$
which was to be proved. $\square $

{\scriptsize \noindent \\
{\sc Justyna Jarczyk\\
Faculty of Mathematics, Computer Science and Econometrics,
University of Zielona G\'ora\\
Szafrana 4a, PL-65-516 Zielona G\'ora, Poland }
\\
e--mail: {\tt j.jarczyk@wmie.uz.zgora.pl}
}

{\scriptsize \noindent \\
{\sc Janusz Matkowski\\
Faculty of Mathematics, Computer Science and Econometrics, 
University of Zielona G\'ora\\
Szafrana 4a, PL-65-516 Zielona G\'ora, Poland }
\\
e--mail: {\tt j.matkowski@wmie.uz.zgora.pl}
}


\begin{thebibliography}{99}
\bibitem{A.Clarc} A. Clarkson, \textit{Uniformly convex spaces}, Trans.
Amer. Math. Soc. 40 (1936), 396-414.
\bibitem{A.Gran} A. Granas, J. Dugundji, \textit{Fixed point theory},
Springer Monographs in Mathematics, Springer Verlag, New York - Berlin -
Heidelberg, 2003.
\bibitem{HLP} G.H. Hardy, J.E. Littlewood, G. P\'{o}lya, \textit{Inequalities}, 
University Press, Cambridge, 1952.
\bibitem{JM89} J. Matkowski, \textit{On a characterization of  $L^{p} $-norm}, 
Ann. Polon. Math. 50 (1989), 137-144.
\bibitem{JMPAMS} -, \textit{The converse of the Minkowski inequality theorem
and its generalization,} Proc. Amer. Math. Soc. 109 (1990), 663-675.
\bibitem{JM.Graz} -, $L^{p} $\textit{-like paranorm}, \textit{Selected Topics
in Functional Equations and Iteration Theory}, Proceedings of the the
Austrian-Polish Seminar, Graz, 1991, Grazer Math. Ber. 316 (1992),
103-138.
\bibitem{JM-Abh} -, \textit{On a generalization of Mulholland's inequality},
Abh. Math. Sem. Univ. Hambg. 63 (1993), 97-103.
\bibitem{JM-Indag} -, \textit{The converse theorem for the Minkowski
inequality}, Indag. Math. (N.S.) 15 (2004), 73-83.
\bibitem{JM-JMAA2008} -, \textit{Converse theorem for the Minkowski
inequality}, J. Math. Anal. Appl. 348 (2008), 277-287.
\bibitem{JM-JMAA-2012} -, \textit{Pexider type generalization of the
Minkowski inequality}, J. Math. Anal. Appl. 393 (2012), 298-310.
\bibitem{Mulhol} H. P. Mulholland, \textit{On generalization of Minkowski
inequality in the form of a triangle inequality,} Proc. London Math. Soc. 
51 (1950), 294-307.
\bibitem{Musielak} J. Musielak, \textit{Introduction to functional analysis}
(in Polish), PWN, Warszawa, 1989.
\bibitem{Pasicki1} L. Pasicki, \textit{A basic fixed point theorem}, Bull.
Pol. Acad. Sci. Math. 54 (2006), 85-88.
\bibitem{Pasicki2} L. Pasicki,\textit{\ Bead spaces and fixed point
theorems}, Topology Appl. 156 (2009), 1811-1816.
\bibitem{Pasicki4} L. Pasicki,\textit{\ Uniformly convex spaces, bead
spaces, and equivalence conditions}, Czechoslovak Math. J. 61 (2011), 383-388.
\bibitem{Pasicki3} L. Pasicki,\textit{\ Towards Lim}, Topology Appl. 156 (2011), 470-483.
\bibitem{Wilansky} A. Wilansky, \textit{Functional analysis}, Blaisdell
Publishing Company, New York - Toronto - London, 1964.
\bibitem{Wnuk} W. Wnuk,\textit{\ Orlicz spaces cannot be normed analogously
to } $L^{p} $\textit{\ spaces},  Nederl. Akad. Wetensch. Indag. Math. 46 (1984), 357-359.
\end{thebibliography}
\end{document}